\documentstyle[12pt,amstex]{article}
\def\IC{\hbox{\rm C \kern-.80em \vrule depth 0ex height 1.5ex width .05em \kern.41em}}
\def\IR{\hbox{\rm R \kern-.80em \vrule depth 0ex height 1.5ex width .05em \kern.41em}}
\def\narrow{\hbox{$\,\,\,\not$ \kern-.80em $\Longrightarrow$}}

\makeatletter
\renewcommand{\section}{\@startsection{section}{1}{0mm}{15mm}{5mm}{\bf\raggedright}}
\makeatother

\newtheorem{theorem}{Theorem}[section]

\newtheorem{coro}[theorem]{Corollary}
\newtheorem{lemma}[theorem]{Lemma}

\newtheorem{remark}[theorem]{Remark}
\newtheorem{example}[theorem]{Example}

\begin{document}
\vspace*{-1.5cm}
\begin{center}
{\large LINEAR EQUATIONS OVER CONES AND \\COLLATZ-WIELANDT NUMBERS} \vspace{0.4cm}\\
{Bit-Shun Tam*}
\vskip 0 pt
{Department of Mathematics}

{Tamkang University}

{Tamsui, Taiwan 25137}

{R.O.C.}\vspace{1mm}\\
and\vspace{1mm}\\
{Hans Schneider}
\vskip 0 pt
{Department of Mathematics}

{University of Wisconsin-Madison}

{Madison, Wisconsin 53706}

{U.S.A.}\vspace{1mm}\\
 {11 September 2001}

\end{center}

\vskip 0 pt
\par\parindent=0 pt

Abstract.$\,\,$ Let $K$ be a proper cone in $\IR^n$,
let $A$ be an $n\!\times\! n$ real matrix that satisfies
$AK\subseteq K$,
let $b$ be a given vector of $K$,
and let $\lambda$ be a given positive real number.
The following two linear equations are considered in this paper:
$\,\,$(i)$\,$ $(\lambda I_n-A)x=b$,
$x\in K$,
and $\,\,$(ii)$\,$ $(A-\lambda I_n)x=b$,
$x\in K$.
We obtain several equivalent conditions for the solvability of the first equation.
For the second equation we give an equivalent condition for its solvability in case when $\lambda>\rho_b (A)$,
and we also find a necessary condition when $\lambda=\rho_b (A)$ and also when $\lambda < \rho_b(A)$,
sufficiently close to $\rho_b(A)$,
where $\rho_b (A)$ denotes the local spectral radius of $A$ at $b$.
With $\lambda$ fixed,
we also consider the questions of when the set $(A-\lambda I_n)K \bigcap K$ equals $\{0\}$ or $K$,
and what the face of $K$ generated by the set is.
Then we derive some new results about local spectral radii and Collatz-Wielandt sets (or numbers) associated with a cone-preserving map,
and extend a known characterization of $M$-matrices among $Z$-matrices in terms of alternating sequences.

\par\parindent=0 pt
\vskip 8 pt
\unitlength 1mm
\begin{picture}(150,5)
\put(0,2){\line(1,0){40}}
\end{picture}
\par\parindent= 0pt
\vspace*{-0.5cm}
\leftmargini=0mm
\begin{enumerate}
\item[]
*Research of this author partially supported by the National Science Council of the Republic of China
{\it E-mail addresses}:
bsm01@@mail.tku.edu.tw
(B.S. Tam);
hans@@math.wisc.edu
(H. Schneider)
\end{enumerate}

\noindent
\newpage
\par\parindent=0 pt
\baselineskip 17 pt{

\section{Introduction}
\par\parindent=16 pt

Let $K$ be a proper
(i.e.,
closed,
pointed,
full,
convex)
cone in $\IR^n$.
Let $A$ be an $n\!\times\! n$ real matrix that satisfies $AK\subseteq K$.
Let $b$ be a nonzero vector of $K$,
and let $\lambda$ be a given positive real number.
Also let $I_n$ denote the $n\!\times\! n$ identity matrix.
In this paper we shall study the solvability of the following two linear equations:
\begin{equation}
(\lambda I_n-A)x=b,\,\,\, x\in K.
\end{equation}

and
\begin{equation}
(A-\lambda I_n)x=b,\,\,\, x\in K.
\end{equation}

Equation (1.1) has been treated before by a number of people.
The study began with the work of Carlson [Car] in 1963 for the special case when $K$ equals $\IR^n_+$
(the nonnegative orthant of $\IR^n$)
and $\lambda$ equals $\rho(A)$
(the spectral radius of $A$),
and was followed by Nelson [Nel 2,3],
Friedland and Schneider [F--S],
Victory [Vic 3],
F$\ddot{\rm o}$rster and Nagy [F--N 1],
and Jang and Victory [J--V 1,2,4].
Indeed,
much of the work on equation (1.1) has been done in the infinite dimensional setting,
when $A$ is either a positive eventually compact linear operator on a Banach lattice,
or an eventually compact linear integral operator with a nonnegative kernel on $L^P(\mu)$ with $1\leq p < \infty$,
or a positive linear operator on an ordered Banach space (or a Banach lattice).
Our contribution here is to provide a more complete set of equivalent conditions for solvability,
and to give simpler and more elementary proofs for the finite dimensional case.

The study of equation (1.2) is relatively new.
A treatment of the equation (by graph-theoretic arguments) for the special case when $\lambda=\rho(A)$ and $K=\IR^n_+$ can be found in Tam and Wu [T--W].
As we shall see,
the solvability of equation (1.2) is a more delicate problem.
It depends on whether $\lambda$ is greater than,
equal to,
or less than $\rho_b(A)$.
When $\lambda$ is fixed,
it is clear that the set $(A -\lambda I_n) K \bigcap K$ consists of preciesly all vectors $b \in K$ for which equation (1.2) has a solution.
(A similar remark can also be said for the set $(\lambda I_n -A)K \bigcap K$.)
For $\lambda \ge 0$,
in general,
the set $(A -\lambda I_n)K \bigcap K$ is not a face of $K$
(but the set $(\lambda I_n -A)K \bigcap K$ is).
So it is natural to consider the face of $K$ generated by $(A -\lambda I_n)K \bigcap K$.
Such a result yields a necessary condition for equation (1.2) to be solvable.
In particular,
in the nonnegative matrix case,
it leads to a combinatorial condition.
We also consider the two extreme situations for the set $(A -\lambda I_n)K \bigcap K$,
namely,
when it is equal to $\{0\}$ or $K$.
As a by-product we obtain a new sufficient condition for $(A -\rho(A) I_n)K \bigcap K=\{0\}$,
which is also a necessary condition in case $K$ is polyhedral.
Note that the condition $(A-\rho(A) I_n)K \bigcap K =\{0\}$ can be rewritten as
``$x \ge^K 0$ and $A x \ge^K \rho(A) x$ imply that $Ax=\rho(A) x$".
Its dual condition ``$x \ge^K 0$ and $A x ~^K\!\!\le \rho(A) x$ imply $Ax =\rho(A) x$" and its equivalent conditions
(see [Tam 1,
Theorem 5.1]) were known and have proved to be useful.

Our work will rely on concepts or results obtained in our recent sequence of papers on the spectral theory of a cone-preserving map ([T--W],
[Tam 1] and [T--S 1,
2]).
(For an overview,
see also [Tam 2].)
In particular,
we frequently make use of a result given in [T--S 2,
Theorem 4.7],
which is about the representation of a nonzero vector $x$ of $K$ in terms of the generalized eigenvectors of the cone-preserving map $A$.
We shall refer to the said representation as the local Perron-Schaefer condition on $A$ at $x$.

Based on our knowledge of the solvability of equations (1.1) and (1.2)
and the local Perron-Schaefer conditions on a cone-preserving map $A$,
we are able to obtain some new results about local spectral radii and Collatz-Wielandt sets (or numbers) associated with $A$.
In particular,
we obtain equivalent conditions for $R_A(x)=\rho_x (A)$,
where $0 \ne x \in K$,
and then characterize when $\inf \sum_1 (=\rho(A)) \in \sum_1$.
For the questions of when $r_A(x)=\rho_x (A)$ and when $\sup \Omega_1 \in \Omega_1$,
we give some partial results.
We also extend a known characterization of $M$-matrices among $Z$-matrices in terms of alternating sequences.

It would be of interest to explore to what extent the methods used in this paper can be carried over to the infinite dimensional settings.

This paper is based on a talk entitled ``Linear equations over cones,
Collatz-Wielandt numbers and local Perron-Schaefer conditions",
given by the first author at the Oberwolfach Workshop on ``Nonnegative matrices,
$M$-matrices and their generalizations"
on November 26--December 2,
2000.
An initial version of this work was also presented by him in the talk ``Solutions of linear equations over cones" at the 11th Haifa Matrix Conference on June 21--25,
1999.

\section{Preliminaries}

We shall restrict our attention to finite-dimensional vector spaces and treat linear equations over proper cones.
By a {\it proper cone} in a finite-dimensional real vector space we mean a nonempty subset $K$ which is a convex cone
(i.e. $\alpha K+\beta K\subseteq K$ for all $\alpha, \beta\geq 0$),
is pointed (i.e. $K \cap (-K)=\{ 0\}$),
has nonempty interior and is closed
(relative to the usual topology of the underlying space).

Hereafter we always use $K$ to denote a proper cone in the Euclidean space $\IR^n$,
and use $\pi(K)$ to denote the set of all $n\!\times\!n$ real matrices $A$ that satisfy $AK\subseteq K$.
(Vectors in $\IR^n$ are represented by $n \times 1$ column vectors.)
Elements of $\pi(K)$ are usually referred to as cone-preserving maps
(or positive operators) on $K$.
It is clear that $\pi(\IR^n_+)$ is the set of all $n\!\times\!n$ nonnegative matrices.

A familiarity with convex cones,
cone-preserving maps,
and graph-theoretic properties of nonnegative matrices is assumed.
For references,
see [Bar],
[B--P] and [Sch 3].
For convenience and to fix notation,
we collect below some of the necessary definitions and known results.

Let $\geq^K$ (also$~^{K}\!\!\leq$) denote the partial ordering of $\IR^n$ induced by $K$,
i.e. $x\geq^K y$ (or $y~^K\!\!\!\leq x$) if and only if $x-y\in K$.
A subset $F$ of $K$ is called a {\it face} of $K$ if it is a convex cone and in addition possesses the property that $x\geq^K y \geq^K 0$ and $x\in F$ imply $y\in F$.
For any subset $S$ of $K$,
we denote by $\Phi(S)$ the {\it face of K generated by S,}
that is,
the intersection of all faces of $K$ including $S$.
If $x\in K$,
we write $\Phi(\{x\})$ simply as $\Phi(x)$.

By the {\it dual cone} of $K$,
denoted by $K^*$,
we mean the (proper) cone $\{z \in \IR^n:~ z^T x \geq 0$ for all $x \in K \}$.
It is well-known that for any $n \times n$ real matrix $A$,
$A \in \pi (K)$ if and only if $A^T \in\pi (K^*)$.

The lattice of all faces of $K$
(under inclusion as the partial ordering) is denoted by ${\cal F}(K)$.
By the duality operator of $K$,
denoted by $d_K$,
we mean the mapping from ${\cal F}(K)$ to ${\cal F}(K^*)$ given by $d_K (F) =({\mbox{span }F})^\bot \bigcap K^*$.
If $F$ is a face of $K$,
we call $d_k(F)$ the {\it dual face} of $F$.
We shall use tacitly the elementary properties of faces and of duality operators.
In particular,
the following fact will be used a number of times:
If $S$ and $T$ are mutually orthogonal nonempty subsets of $K$ and $K^*$ respectively,
then $\Phi(T)\subseteq d_K (\Phi (S))$ and $\Phi(S) \subseteq d_{K^*} (\Phi(T))$.

Let ${\cal M}_n (\IC)$ denote the space of all $n\!\times\!n$ complex matrices,
and let $A \in {\cal M}_n(\IC)$.
The range space,
nullspace and the spectral radius of $A$ are denoted respectively by ${\cal R}(A),{\cal N}(A)$ and $\rho(A)$.
Eigenvalues of $A$ with modulus $\rho (A)$ are said to compose the {\it peripheral spectrum of $A$}.
For any eigenvalue $\lambda$ of $A$,
we use $\nu_{\lambda}(A)$ to denote the {\it index of $\lambda$ as an eigenvalue of $A$},
i.e.,
the smallest integer $k$ such that ${\cal N}((A-\lambda I_n)^k)={\cal N}((A-\lambda I_n)^{k+1})$.
For any vector $x\in \IC^n$,
by the {\it cyclic space relative to A generated by x},
denoted by $W_x$,
we mean the linear subspace span$\{x, Ax, A^2 x,\ldots\}$.

We also use ${\cal M}_n(\IR)$ to denote the space of all $n\!\times\!n$ real matrices.
The above concepts and notation will also apply to the real case.
Sometimes we treat an $n \times n$ real matrix $A$ as a complex matrix;
in other words,
we identify $A$
(as a linear operator)
with its complex extension acting in the complexification $\IC^n$ of $\IR^n$.
So for $A\in {\cal M}_n (\IR)$,
the symbol ${\cal N}((\rho (A) I_n -A)^n)$ (also $W_x$ with $x \in \IR^n$)
can represent a real subspace of $\IR^n$ or a complex subspace of $\IC^n$,
as understood from the context.

We need the concept of the {\it local spectral radius of A at x},
which is denoted by $\rho _x (A)$.
If $x$ is the zero vector,
take $\rho_x (A)$ to be 0.
Otherwise,
define $\rho _x (A)$ in one of the following equivalent ways
(see [T--W,
Theorem 2.3]):

(i) $\rho_x (A) =\lim \sup_{m\rightarrow\infty} \|A^m x\|^{1/m}$,
where $\| \cdot \|$ is any norm of $\IC^n$.

(ii) $\rho _x (A) =\rho (A|_{W_x})$,
where $A|_{W_x}$ denotes the restriction of $A$ to the invariant subspace $W_x$.

(iii) Write $x$ uniquely as a sum of generalized eigenvectors of $A$,
say,
$x=x_1 +\cdots+x_k$,
where $k \geq 1$ and $x_1,\ldots,x_k$ are generalized eigenvectors of $A$ corresponding to distinct eigenvalues $\lambda_1,\ldots,\lambda_k$.
Then define $\rho_x (A)$ to be $\max _{1\leq i \leq k} |\lambda_i|$.

It is worth noting that,
in the first definition of local spectral radius,
we can replace ``lim sup" by ``lim".
We shall offer a proof for this assertion in Appendix A to this paper.

We also need the concept of the order of a vector relative to a square matrix as introduced in [T--S 2].
Let $A \in {\cal M}_n (\IC)$.
If $x$ is a generalized eigenvector of $A$ corresponding to the eigenvalue $\lambda$,
then by the order of $x$ we mean,
as usual,
the least positive integer $p$ such that $(A-\lambda I_n)^p x=0$.
If $x$ is a nonzero vector of $\IC^n$,
then by the {\it order of x relative to A},
denoted by ${\rm ord}_A (x)$,
we mean the maximum of the orders of the generalized eigenvectors,
each corresponding to an eigenvalue of modulus $\rho_x (A)$,
that appear in the representation of $x$ as a sum of generalized eigenvectors of $A$.
(If $x$ is the zero vector,
we set ${\rm ord}_A (x)$ to be 0.)

It is convenient to introduce the concept of spectral pair here.
Following [T--S 2],
we denote the ordered pair $(\rho_x (A), {\rm ord}_A(x))$ by sp$_A (x)$ and refer to it as the {\it spectral pair of $x$ relative to $A$}.
We also denote by $\preceq$ the lexicographic ordering between ordered pairs of real numbers given by:
\[
( \xi_1, \xi_2) \preceq (\eta_1, \eta_2)\mbox{ if either }\xi_1 < \eta_1\mbox{ or }\xi_1 =\eta_1\mbox{ and }\xi_2 \le \eta_2.
\]
In [T--S 2,
Theorem 4.7] it is shown that if $A \in \pi (K)$,
then for any face $F$ of $K$,
the spectral pair sp$_A(x)$ is independent of the choice of $x$ from relint $F$.
This common value is denoted by sp$_A(F)$ and is called the {\it spectral pair of $F$ relative to $A$}.
The concept of spectral pair of faces (or vectors) plays an important role in the work of [T--S 2].

Let $A\in {\cal M}_n (\IR)$.
It is known that a necessary and sufficient condition for the existence of a proper cone $K$ of $\IR^n$ such that $A\in \pi (K)$ is that the following set of conditions is satisfied:

(a)
$\rho(A)$ is an eigenvalue of $A$.

(b)
If $\lambda$ is an eigenvalue in the peripheral spectrum of $A$,
then $\nu_{\lambda} (A) \leq \nu_{\rho(A)} (A)$.\\
The above set of conditions is now referred to as the {\it Perron-Schaefer condition}
(see [Sch 2,
the paragraph following Theorem 1.1],
and also [T--S 1,
Section 7] for our recent work involving the condition).

According to [T--S 2,
Theorem 4.7],
if $A\in \pi (K)$,
then for any $0 \neq x \in K$,
the following condition is always satisfied:

{\it
There is a generalized eigenvector $y$ of $A$ corresponding to $\rho_x (A)$ that appears as a term in the representation of $x$ as a sum of generalized eigenvectors of $A$.
Furthermore,
we have ord$_A (x)={\rm ord}_A (y)$.}

By analogy,
we shall refer to the preceding condition as the {\it local Perron-Schaefer condition on $A$ at $x$}.
In a forthcoming paper [Tam 3] we shall show that $A$ satisfies the local Perron-Schaefer condition at $x$ if and only if $A|_{W_x}$ satisfies the Perron-Schaefer condition.
Still another equivalent condition is that,
the closure of the convex cone generated by $A^ix$ for $i=0, 1, \ldots$ is pointed.
Based on the equivalence of these conditions and a result of similar kind,
in [Tam 3] we also rederive certain intrinsic Perron-Frobenius theorems obtained by Schneider in [Sch 2].

If $A \in \pi (K)$ and $x \in K$ is an eigenvector
(respectively,
generalized eigenvector),
then $x$ is called a {\it distinguished} eigenvector
(respectively,
{\it distinguished generalized eigenvector})
{\it of A for K},
and the corresponding eigenvalue is known as a {\it distinguished eigenvalue of A for K}.
When there is no danger of confusion,
we simply use the terms distinguished eigenvector
(respectively,
distinguished generalized eigenvector)
and distinguished eigenvalue (of $A$).
It is known that a real number $\lambda$ is a distinguished eigenvalue of $A$ if and only if $\lambda = \rho_x (A)$ for some $0 \neq x \in K$
(see [T--W, Theorem 2.4]).

Let $A \in \pi (K)$.
A face $F$ of $K$ is said to be {\it A-invariant} if $AF \subseteq F$.
The following result is proved in [T--S 2,
Lemma 2.1,
Corollary 4.5 and Theorem 4.9 (ii) (a)]:

\begin{lemma}
Let $A\in\pi(K)$ and let $x\in K$.
Also let ${\hat x}={(I+A)}^{n-1}x$.
Then $\Phi({\hat x})$ is the smallest $A$-invariant face of $K$ containing $x$,
and $W_x={\rm span} \Phi(\hat{x})$.
Furthermore,
$\rho_x (A)=\rho_{\hat{x}}(A) =\rho( A|_{{\rm span}\Phi(\hat{x}) } )$ and
${\rm ord}_A(x)={\rm ord}_A(\hat{x})$.
\end{lemma}

Hereafter we use $P$ to denote an $n\!\times\!n$ nonnegative matrix.
The set $\{1,2,\cdots,n\}$ is denoted by $\langle n \rangle$.
For any nonempty subsets $I,J$ of $\langle n \rangle$,
we use $P_{IJ}$ to denote the submatrix of $P$ with rows indexed by $I$ and columns indexed by $J$.
We follow the standard usage of the concepts of classes of $P$ and of accessibility relation and denote the classes by Greek letters $\alpha, \, \beta$,
etc.
(see [Rot] or [Sch 3]).
The accessibility relation is usually defined between the classes of $P$.
For convenience,
we also allow the relation be defined in the natural way between the elements of $\langle n\rangle$,
between the nonempty subsets of $\langle n \rangle$,
and between the elements of $\langle n \rangle$ and the nonempty subsets of $\langle n\rangle$.
For instance,
if $i \in \langle n \rangle$ and $\emptyset \neq J \subseteq \langle n \rangle$,
we say $i$ {\it has access} to $J$ if there is a path in the directed graph of $P$ from the vertex $i$ to some vertex in $J$.
If $\alpha, ~\beta$ are classes of a nonnegative matrix $P$,
we write $\alpha >= \beta$ if $\alpha$ has access to $\beta$.
We also write $\alpha >\!\!- \,\beta$ if $\alpha >= \beta$ but $\alpha \ne \beta$.

We also need the concept of an initial subset for $P$ as introduced in [T--S 2].
A subset $I$ of $\langle n \rangle$ is called an {\it initial subset for} $P$ if either $I$ is empty,
or $I$ is nonempty and $P_{I^{\prime}I} =0$,
where $I^{\prime}=\langle n \rangle \backslash \, I$;
equivalently,
for every $j \in \langle n \rangle$,
$I$ contains $j$ whenever $j$ has access to $I$.
It is not difficult to show that a nonempty subset $I$ of $\langle n \rangle$ is an initial subset for $P$ if and only if $I$ is the union of an {\it initial collection of classes} of $P$,
where a nonempty collection of classes of $P$ is said to be initial if whenever it contains a class $\alpha$,
it also contains every class that has access to $\alpha$.

We follow the usual definitions of basic
(initial,
final,
distinguished) class of a nonnegative matrix $P$.
A class $\alpha$ is said to be {\it semi-distinguished} if $\rho(P_{\beta \beta}) \le \rho(P_{\alpha \alpha})$ for any class $\beta >= \alpha$.
For convenience,
we say a class $\alpha$ is associated with $\lambda$ if $\rho (P_{\alpha \alpha})=\lambda$.
If $\cal L$ is a collection of classes of $P$,
then we also say a class $\alpha \in {\cal L}$ is {\it final in}
(respectively,
{\it initial in}) $\cal L$ if $\alpha$ has no access to
(respectively,
access from) any other class in $\cal L$.

It is well-known that every face of $\IR^n_+$ is of the form
\[
F_I = \{ x \in \IR^n_+ : {\rm supp} (x) \subseteq I\},
\]
where $I \subseteq \langle n \rangle$,
and ${\rm supp} (x)$ is the {\it support of x},
i.e. the set $\{ i \in \langle n \rangle : \xi_i \neq 0 \}$ for $x=(\xi_1,\cdots,\xi_n)^T$.

We need the following result which is proved in [T--S 2,
Theorem 3.1]:

\begin{theorem}
Let P be an $n \! \times\! n$ nonnegative matrix.
Denote by ${\cal F}_P$ the lattice of all P-invariant faces of $\IR^n_{+}$ and by $\cal I$ the lattice of all initial subsets for P,
both under inclusion as the partial ordering.
Then the association $I \longmapsto F_I$ induces an isomorphism from the lattice $\cal I$ onto the lattice ${\cal F}_P$.
\end{theorem}

For any $A \in \pi (K)$,
the following sets are called the {\it Collatz-Wielandt sets associated with $A$}:
\begin{eqnarray*}
\Omega (A) &=& \{ \omega \geq 0: \exists x \in K \backslash \{ 0 \}, Ax \geq^{K} \omega x \}, \\
\Omega_1 (A) &=& \{ \omega \geq 0: \exists x \in {\rm int}\, K, Ax \geq^{K} \omega x \},\\
{\textstyle\sum} (A) &=& \{\sigma \geq 0 : \exists x \in K \backslash \{ 0 \}, Ax ~^K\!\!\!\leq \sigma x\},\\
{\textstyle\sum_1} (A) &=& \{ \sigma \geq 0 : \exists x \in {\rm int}\, K, Ax ~^K\!\!\!\leq \sigma x \}.
\end{eqnarray*}
When there is no danger of confusion,
we write simply $\Omega, \Omega_1, \sum$ and $\sum_1$.
If $x \in K$,
then the {\it lower and upper Collatz-Wielandt numbers of $x$ with respect to A} are defined by
\begin{eqnarray*}
r_A (x) &=& \sup \{ \omega \geq 0: Ax \geq^K \omega x \}, \\
R_A (x) &=& \inf \{ \sigma \geq 0: Ax ~^K\!\!\!\leq \sigma x \},
\end{eqnarray*}
where we write $R_A (x)=\infty$ if no $\sigma$ exists such that $Ax ~^K\!\!\!\leq \sigma x$.
It is clear that when equation (1.1) is solvable,
we have $\lambda \in \sum$ and $\lambda \geq R_A (x)$ for any solution $x$.
We refer our reader to [T--W] for results on the Collatz-Wielandt sets or numbers.
(See also [F--N 2],
[Fri] and [Mar] for results in the infinite dimensional settings.)

The following known result will be used tacitly:

\begin{remark}
{\rm
Let $A \in {\cal M}_n (\IC)$.
For any eigenvalue $\lambda$ of $A$,
the orthogonal complement in $\IC^n$ of the generalized eigenspace of $A$ corresponding to $\lambda$ is equal to the direct sum of all generalized eigenspaces of $A^*$ corresponding to eigenvalues other than $\bar{\lambda}$.}
\end{remark}

To prove this,
use the following two facts,
valid for any $A \in {\cal M}_n (\IC)$:
(i) ${\cal N}(A)^\bot ={\cal R}(A^*)$;
and (ii) $\IC^n$ is the direct sum of all generalized eigenspaces of $A$.

\section{The equation {\boldmath $(\lambda I_n -A)x=b, x \in K$}.}

\begin{theorem}
Let $A \in \pi (K)$,
let $0 \neq b \in K$,
and let $\lambda$ be a positive real number.
The following conditions are equivalent$\,:$
\vspace*{-2mm}
\leftmargini=6mm
\begin{enumerate} \itemsep=-1pt
\item[] {\rm (a)}
There exists a vector $x\in K$ such that $(\lambda I_n -A)x=b$.
\item[] {\rm (b)}
$\rho_b (A) < \lambda.$
\item[] {\rm (c)}
$\lim\limits_{m\to \infty} \sum\limits^m\limits_{j=0} \lambda^{-j} A^j b$ exists.
\item[] {\rm (d)}
$\lim\limits_{m\to \infty} ({\lambda}^{-1} A)^m b=0$.
\item[] {\rm (e)}
$\langle z,b \rangle =0$ for each generalized eigenvector z of $A^T$ corresponding to an eigenvalue with modulus greater than or equal to $\lambda$.
\item[] {\rm (f)}
$\langle z,b \rangle =0$ for each generalized eigenvector z of $A^T$ corresponding to a distinguished eigenvalue of A for K which is greater than or equal to $\lambda$.
\end{enumerate}
When the equivalent conditions are satisfied,
the vector $x^0=\sum^{\infty}_{j=0} {\lambda}^{-j-1} A^j b$ is a solution of the equation $(\lambda I_n -A)x=b,x \in K$.
Furthermore,
if $\lambda$ is a distinguished eigenvalue of $A$,
then the solution set of the equation consists of precisely all vectors of the form $x^0 + u$,
where u is either the zero vector or is a distinguished eigenvector of A corresponding to $\lambda ;$
otherwise,
$x^0$ is the unique solution of the equation.
\end{theorem}

\noindent
{\it Proof.}

(b)$\Longrightarrow$(c):
Since $\rho(A|_{W_b})= \rho_b (A) < \lambda, (\lambda I_n -A|_{W_b} )^{-1}$ exists and is given by:
\[
(\lambda I_n -A|_{W_b} )^{-1}= \lambda^{-1} \lim_{j \to \infty} \sum^m_{j=0} (\lambda^{-1} A|_{W_x})^j,
\]
hence (c) follows.

(c)$\Longrightarrow$(d):
Obvious.

(d)$\Longrightarrow$(b):
Condition (d) clearly implies that
\[
\lim_{m\to \infty} (\lambda^{-1} A)^m (A^i b)=0 \mbox{~~for~~} i=0,1,2,\cdots \, ,
\]
and hence
\[
\lim_{m\to \infty} (\lambda^{-1} A)^m y=0 \mbox{~~for every~~} y \in W_b .
\]
It follows that ~$\lim_{m\to \infty} (\lambda^{-1} A|_{W_b})^m =0,\mbox{~and hence~~} \rho_b (A) < \lambda.$

(b)$\Longrightarrow$(e):
Since $\rho_b (A)<\lambda$,
$b$ is contained in the direct sum of all generalized eigenspaces of $A$ corresponding to eigenvalues with moduli less than $\lambda$.
It follows that for any generalized eigenvector $z$ of $A^T$ corresponding to an eigenvalue with modulus greater than or equal to $\lambda$,
we have $\langle z,b \rangle=0$.

(e)$\Longrightarrow$(b):
Suppose that $\rho_b(A) \geq \lambda$.
Let $\mu$ be an eigenvalue of $A|_{W_b}$ such that $|\mu|=\rho_b(A)$.
Since $A$ is a real matrix,
$\bar{\mu}$ is also an eigenvalue of $A$ and hence of $A^T$.
According to condition (e),
$b$ belongs to,
and hence $W_b$ is included in,
the $A$-invariant subspace $({\cal N}((\bar{\mu} I_n -A^T)^n))^{\bot}$ of $\IC^n$,
which is the same as the direct sum of all generalized eigenspaces of $A$ corresponding to eigenvalues other than $\mu$.
Clearly $\mu$ is not an eigenvalue of the restriction of $A$ to the latter subspace.
On the other hand,
by our choice,
$\mu$ is an eigenvalue of $A|_{W_b}$.
So we arrive at a contradiction.

We have just shown that for any $A \in M_n ( \IR )$ and $0 \neq b \in \IR^n$,
conditions (b), (c), (d) and (e) are equivalent.

(c)$\Longrightarrow$(a):
Obvious:
put $x=\sum^{\infty}_{k=0} \lambda^{-(k+1)} A^k b$.

(a)$\Longrightarrow$(c):
Suppose that there exists a vector $x\in K$ such that $(\lambda I_n -A)x=b$.
Then we have,
$x=\lambda^{-1} b+\lambda^{-1} Ax \geq^K \lambda^{-1} b$,
where the inequality follows from the assumptions that $A \in \pi (K)$ and $x\in K$.
Multiplying both sides of the inequality by $\lambda^{-1}A$,
we obtain $\lambda^{-1} Ax \geq^K \lambda^{-2} Ab$,
and hence $x \geq^K \lambda^{-1} b +\lambda^{-2} Ab$.
Proceeding inductively,
we show that,
for all positive integers $m$,
we have $x\geq^K y_m$,
where we denote by $y_m$ the vector $\sum^m_{j=0} {\lambda}^{-j-1} A^j b$.
(Actually,
as noted in [F--N 1, Remark 9],
we have $x=\lambda^{-m} A^m x+y_m$ for all positive integers $m$.)
So we have $0 ~^K\!\!\!\leq y_1 ~^K\!\!\!\leq y_2 ~^K\!\!\!\leq \cdots ~^K\!\!\!\leq x$,
and from this we are going to deduce that $\lim y_n$ exists.
(This is,
undoubtedly,
a known fact.
See,
for instance,
[G--L,
p.355,
Problem 20].
For completeness,
we supply a proof here.)

Choose a norm $\| \cdot \|$ of $\IR^n$ which is monotonic with respect to $K$
(see [B--P, p.6]).
Then we have $\| y_m \| \leq \|x \|$ for all $m$.
So $(y_m)$ is a bounded sequence of $\IR^n$,
and we can choose a convergent subsequence,
say,
$(y_{k_j})$ with limit $y$.
Then for any positive integer $m$,
we have $y_{k_j} \geq^K y_m$ whenever $k_j \geq m$,
and by letting $j \to \infty$,
we obtain $y \geq^K y_m$.
Indeed,
we have $\lim y_m =y$;
this is because for any $m \geq k_j$,
we have,
$0 ~^K\!\!\!\leq y-y_m ~^K\!\!\!\leq y- y_{k_j}$,
and hence $\|y - y_m \| \leq \| y - y_{k_j}\|$.
So the desired limit exists.

(e)$\Longrightarrow$(f):
Obvious.

(f)$\Longrightarrow$(b):
This follows from the fact that $\rho_b (A)$ is a distinguished eigenvalue of $A$ for $K$ as $0\neq b \in K$
(see [T--W,
Theorem 2.4(ii)]) and also that $b \not\in ({\cal N}((\rho_b (A) I_n -A^T)^n))^{\bot}$
(cf. the argument used in the proof of (e) $\Rightarrow$ (b)).

This proves the equivalences of conditions (a)--(f).

{\it Last Part.}
>From the above proof of (a)$\Longrightarrow$(c),
it is clear that,
when the equivalent conditions are satisfied,
the vector $x^0=\sum^{\infty}_{j=0} {\lambda}^{-j-1} A^j b$ is a solution of equation (1.1).
Moreover,
the proof also shows that if $x\in K$ satisfies (1.1),
then $x \geq^K y_m$ for each positive integer $m$,
where $y_m=\sum^m_{j=0} {\lambda}^{-j-1} A^j b$.
But $\lim\limits_{m\to \infty} y_m =x^0$,
hence $x-x^0=u$,
where $u\in K$.
Since $x$ and $x^0$ both satisfy (1.1),
it is clear that $u$ is either the zero vector or is an eigenvector of $A$ corresponding to $\lambda$.
Hence,
our assertion follows.\hfill $\Box$\\

In the proof of Theorem 3.1,
instead of proving the implication (a) $\Rightarrow$ (c),
we can also proceed by establishing the implication (a) $\Rightarrow$ (b).
(Then we prove the last part of our theorem by using the argument given in [F--N 1,
Remark 9].)
We have found two interesting proofs for the implication (a) $\Rightarrow$ (b).
We include them in Appendix B to this paper.

The following alternative proof of Theorem 3.1,
(f)$\Longrightarrow$(b) that makes use of the local Perron-Schaefer condition is also of interest:

Suppose condition (b) does not hold,
i.e.,
$\rho_b(A) \ge \lambda$.
Let $b=b_1+\cdots+b_k$ be the decomposition of $b$ in terms of generalized eigenvectors of $A$.
By the local Perron-Schaefer condition on $A$ at $b$,
we may assume that the generalized eigenvector $b_1$ corrresponds to the eigenvalue $\rho_b(A)$.
Choose a generalized eigenvector $z$ of $A^T$ corresponding to $\rho_b(A)$ such that $\langle z, b_1 \rangle \ne 0$.
Then $\langle z, b \rangle=\langle z, b_1 \rangle \ne 0$.
But $\rho_b (A)$ is a distinguished eigenvalue of $A$ for $K$,
hence condition (f) does not hold.\\

The first half of the proof of Theorem 3.1 actually shows the following:

\begin{remark}
{\rm
For any $A \in {\cal M}_n (\IC)$,
any $b \in \IC^n$,
and any positive real number $\lambda$,
conditions (b), (c), (d) of Theorem 3.1 and the following condition (e)$^{\prime}$ are equivalent:

{(e)$^{\prime}$}
$\langle z,b \rangle=0$ for each generalized eigenvector $z$ of $A^*$ corresponding to

an eigenvalue with modulus greater than or equal to $\lambda$.\\
Moreover,
the following condition is always implied by the above equivalent conditions,
but is not equivalent to them:

{(a)$^{\prime}$}
There exists a vector $x \in \IC^n$ such that $(\lambda I_n -A)x =b$.}
\end{remark}

Since condition (a)$^{\prime}$ is satisfied whenever $\lambda$ is not an eigenvalue of $A$,
clearly condition (a)$^{\prime}$ does not imply condition (b) of Theorem 3.1.
In fact,
even if $\lambda$ is an eigenvalue of $A$,
the implication (a)$^{\prime} \Longrightarrow$(b) still does not hold.
As a counter-example,
consider $A={\rm diag}(0,1 ,2), \lambda=1$ and $b=(0,0,1)^T$.
In this case,
condition (a)$^{\prime}$ is satisfied,
but we have $\lambda=1<2=\rho_b (A)$.\\

If $z\in \IC^n$, we write $|z|$ to mean the nonnegative vector whose components are the moduli of the corresponding components
of $z$.

\begin{coro}
Let P be an $n\!\times\!n$ nonnegative matrix,
let $b\in\IR^n_+$,
and let $\lambda$ be a positive real number.
To the list of equivalent conditions of Theorem 3.1
$($but with A and K replaced respectively by P and $\IR^n_+ )$
we can add the following$\,:$
\vspace*{-2mm}
\leftmargini=6mm
\begin{enumerate} \itemsep=-1pt
\item[] {\rm (g)}
For any class $\alpha$ of P having access to ${\rm supp}(b)$,
$\rho (P_{\alpha \alpha})< \lambda$.
\item[] {\rm (h)}
For each distinguished class $\alpha$ of $P$ for which $\rho(P_{\alpha \alpha}) \ge \lambda$,
we have $b_\beta =0$ whenever $\beta$ is a class that has access from $\alpha$.
\item[] {\rm (i)}
$\langle |z|, b \rangle =0$ for each generalized eigenvector z of $P^T$ corresponding to an eigenvalue with modulus greater than or equal to $\lambda$.
\item[] {\rm (j)}
$\langle |z|, b \rangle =0$ for each generalized eigenvector z of $P^T$ corresponding to a distinguished eigenvalue of P for $\IR^n_+$ which is greater than or equal to $\lambda$.
\end{enumerate}
\vspace*{-2mm}
When the equivalent conditions are  satisfied,
the vector $x^0=\sum\limits^{\infty}\limits_{j=0} \lambda^{-j-1} P^j b$ is a solution for the given equation and is also the unique solution with the property that its support is included in $($in fact,
equal to$)$ the union of all classes of P having access to ${\rm supp}(b)$.
In this case,
if $\lambda$ is not a distinguished eigenvalue of P, then $x^0$ is the unique solution,
and if $\lambda$ is a distinguished eigenvalue,
then the solutions of the equation are precisely all the vectors of the form $x^0+u$,
where u is the zero vector or is a distinguished eigenvector of P corresponding to $\lambda$.
\end{coro}

\noindent
{\it Proof.}
Note that condition (g) can be rewritten as:
For each class $\alpha$ of $P$ for which $\rho(P_{\alpha \alpha}) \ge \lambda$,
we have $b_\beta =0$ whenever $\beta$ is a class that has access from $\alpha$.
So,
clearly we have (g)$\Longrightarrow$(h).
To show the reverse implication,
let $\cal L$ denote the collection of all classes $\alpha$ for which $\rho(P_{\alpha \alpha})\ge \lambda$.
Consider any $\alpha \in {\cal L}$.
Let $\gamma$ be a class initial in $\cal L$ such that $\gamma >= \alpha$.
Then for any class $\delta >\!\!-\, \gamma$,
we have $\delta \notin {\cal L}$ and hence $\rho(P_{\delta \delta})< \lambda$.
Thus,
$\gamma$ is a distinguished class of $P$.
If $\beta$ is a class such that $\alpha >= \beta$,
then we also have $\gamma >= \beta$ and by condition (h) it follows that we have $b_\beta =0$.
This establishes (h)$\Longrightarrow$(g).

It is clear that condition (g) is equivalent to condition (b) once we prove the following assertion:
\[
\rho_b (P)=\max \{ \rho (P_{\alpha \alpha}) :
\alpha \mbox{ has access to supp}(b) \}.
\]
According to Lemma 2.1,
$\rho_b (P)$ is equal to the spectral radius of the restriction of $P$ to the linear span of the
smallest $P$-invariant face of $\IR^n_+$ containing $b$.
Let $I$ denote the smallest initial subset for $P$ including ${\rm supp}(b)$,
i.e.,
the union of all classes of $P$ having access to ${\rm supp}(b)$.
By Theorem 2.2 the smallest $P$-invariant face containing $b$ is clearly equal to $F_I$.
So we have
\begin{eqnarray*}
\rho_b (P) &=& \rho (P|_{{\rm span}\, F_I} ) \\
           &=& \max \{ \rho (P_{\alpha \alpha}) : \alpha \mbox{ has access to supp}(b) \}.
\end{eqnarray*}
This proves our assertion and hence the equivalence of conditions (g) and (b).

Note that,
since $b$ is a nonnegative vector,
condition (j) amounts to saying that ${\rm supp}(b)\bigcap {\rm supp}(z)=\emptyset$,
for any vector $z$ with the property described in (j).
So it is clear that we have the implication (j)$\Longrightarrow$(f).

Clearly we also have the implication (i)$\Longrightarrow$(j).
It remains to show (g)$\Longrightarrow$(i).
Let $\cal C$ denote the collection of all classes of $P$ that have access from some class $\alpha$ for which $\rho(P_{\alpha \alpha}) \geq \lambda$,
and let $J$ denote the union of all classes in $\cal C$.
Let $J^{\prime}$ denote $\langle n \rangle \backslash J$.
It is clear that $P_{J J^{\prime}}=0$.
Also it is not difficult to see that $\rho (P_{J^{\prime} J^{\prime}}) < \lambda$.
By a permutation similarity,
we may assume that
\[
P=\left[
\begin{array}{cc}
P_{JJ} & 0 \\
P_{J^{\prime} J} & P_{J^{\prime} J^{\prime}}
\end{array} \right].
\]
Consider any generalized eigenvector $z$ of $P^T$ corresponding to an eigenvalue $\mu$ with $| \mu | \geq \lambda$.
Partition $z$ as
$\left[
\begin{array}{c}
u \\ v
\end{array} \right]$ in conform with the above partitioning of $P$.
By definition there exists a positive integer $k$ such that $(P^T - \mu I_n)^k z=0$.
A little calculation shows that $((P_{J^{\prime} J^{\prime}})^T -\mu I_n)^k v =0$.
Since $\mu$ cannot be an eigenvalue of $(P^{J^{\prime} J^{\prime}})^T$,
it follows that $v=0$.
This shows that supp$(z) \subseteq J$.
On the other hand,
by condition (g),
supp$(b) \bigcap J =\emptyset$.
Hence we have $\langle |z|, b \rangle =0$.

{\it Last Part.}
In view of the last part of Theorem 3.1,
it suffices to show that $x^0$ is the unique solution of the given equation with the property that supp$(x^0)$ equals the union of all classes of $P$ that have access to supp$(b)$.
As in the beginning part of our proof,
let $I$ denote the union of all classes of $P$ that have access to supp$(b)$.
As shown before,
we have $\rho (P|_{{\rm span}\, F_I} )=\rho_b (P) <\lambda$,
i.e.$\,\lambda I_n -P|_{{\rm span}\, F_I}$ is nonsingular.
Hence,
the given equation admits a unique solution in $F_I$.
By the definitions of $x^0$ and $I$,
it is also clear that supp$(x^0)=I$.
Thus,
our assertion follows.\hfill $\Box$
\\

Some historical remarks in concern with the conditions of Theorem 3.1 and Corollary 3.3 are in order.

\begin{remark}
{\rm
It was Carlson [Car] who first considered the solvability of equation (1.1) for the special case when $K$ equals $\IR^n_+$ and $\lambda =\rho (A)$
(i.e. when $\lambda I_n -A$ is a singular $M$-matrix).
His work was motivated by the results and technique developed by Schneider [Sch 1].
He gave one equivalent condition,
which in modern language is condition (g) of Corollary 3.3.
He also mentioned that when the equation is solvable there is a unique solution such that its support is precisely the union of all classes having access to supp$(b)$.
Motivated by the needs in some fields of applied mathematics
(for instance,
radiative transfer,
linear kinetic theory;
see [Nel 1]),
but apparently unaware of Carlson's work,
Nelson [Nel 2] also considered the solvability of equation (1.1) in the setting when $A$ is a nonnegative eventually compact linear operator and the underlying space is a real vector lattice which is a Banach space under a semi-monotonic norm
(which is more general than a Banach lattice.)
In the Theorem of the paper,
he gave one equivalent condition which corresponds to condition (j) of Corollary 3.3.
In his proof he actually showed that condition (c) of Theorem 3.1 is also another equivalent condition.
Indeed,
Nelson noted that when the equation is solvable,
$x^0$ (of Theorem 3.1)
is one solution and this solution is majorized by every other solution of the equation.
The local spectral radius $\rho_b (A)$ was actually involved in his proof,
but he didn't use the terminology.
His proof is function-theoretic,
making use of the resolvent of $A$ and also certain kind of Pringsheim's theorem.
In the Corollary,
Nelson treated the special case when the underlying space is $L^P (\mu)$,
$1 \leq p < \infty$,
where the measure $\mu$ is totally $\sigma$-finite if $p=1$.
In the nonnegative matrix case,
his condition reduces to condition (h) of Corollary 3.3.
The investigation was continued in the subsequent paper [Nel 3].
Friedland and Schneider [F--S,
Theorem 7.1] also considered equation (1.1) for the special case when $K=\IR^n_+$ and $\lambda =\rho (A)$,
and extended Carlson's result by adding conditions (c),
(d) of Theorem 3.1.
Their proof is based on an analytic result,
Theorem 5.10 of their paper,
about the growth of $A^m (I_n +A+\cdots+A^{q-1}),m=1,2,\cdots$,
where $A \geq 0$ with $\rho (A)=1$ and $q$ is a certain positive integer.
Victory [Vic 3] considered equation (1.1) for the nonnegative matrix case in which $\lambda$ need not be $\rho (A)$,
and established the equivalences of conditions corresponding to those in [F--S,
Theorem 7.1].
He also gave an additional equivalent condition,
which is condition (j) of Corollary 3.3.
Condition (f) of Theorem 3.1 also appears as a remark at the end of Section 2 of his paper.
The proof given in [Vic 3] for his main Theorem is somewhat involved.
In particular,
the proof of the equivalence of conditions (g) and (j) of Corollary 3.3 relies on a knowledge of the support structure of the generalized eigenvectors of a nonnegative matrix.
In [H--S 2,
Theorem 3.11],
Hershkowitz and Schneider also treated equation (1.1) for the nonnegative matrix case.
They formulated their results in terms of Z-matrices and provide a proof
(in fact,
two for the nontrivial direction)
for the equivalence of condition (g) of Corollary 3.3 and condition (a) of Theorem 3.1
(with $K$ and $A$ replaced by $\IR^n_+$ and $P$ respectively).
They also investigated the case when nonnegativity of the solution is not required.
(See Section 4 of their paper.)
About the same time,
F$\ddot{\rm o}$rster and Nagy [F--N 1],
using the iterated local resolvent as a main tool,
also treated equation (1.1) in the setting of a nonnegative linear continuous operator $T$ on an ordered Banach space $E$,
under the assumption that the positive cone of $E$ is normal,
or that the local spectral radius $\rho_x (T)$ is a pole of the local resolvent function $x_T (\cdot)$
(given by $x_T (\mu)=\sum^{\infty}_{j=0} \mu^{-j-1} T^j x$ for $|\mu|>\rho_x (T)$),
or the $\lambda$ of equation (1.1) is greater than the radius of the essential spectrum of $T$
(see [F--N 1,
Theorems 6,
10 and 12]).
In their Theorem 12,
they gave equivalent conditions for the existence of a solution to equation (1.1),
which correspond to conditions (b) and (e) of our Theorem 3.1.
(The condition (iii) in their theorem,
which does not involve $\lambda$,
is a tautology and should be deleted.)
In the early nineties,
the main theorem in [Vic 3] was extended by Jang and Victory first to the setting of an eventually compact linear integral operator with nonnegative kernel,
mapping the Lebesgue space $L^p (\Omega ,\Sigma , \mu),~1 \leq p < \infty$,
into itself,
where $(\Omega ,\Sigma , \mu)$ denotes a $\sigma$-finite measure space
(see [J--V 1,
Theorem 3.6] or [J--V 2,
Theorem IV.1]),
and then to the setting of a nonnegative eventually compact reducible linear operator which maps a Banach lattice having order continuous norm into itself
(see [J--V 4,
Theorem 4.1]).
Their work depends on a knowledge of the structure of the generalized eigenspace corresponding to the spectral radius of a nonnegative reducible eventually compact linear operator defined on a Banach lattice with order continuous norm
(in particular,
a Lebesgue space $L^p (\Omega ,\Sigma , \mu)$).
(The most difficult part in their proof involves a condition in their main result which corresponds to condition (f) of our Theorem 3.1.)
Their investigation is made possible by a decomposition of the underlying operator in a form generalizing the Frobenius normal form of a nonnegative reducible matrix.
(Such decomposition was initiated by Nelson.
For the details,
see [Nel 3],
[Vic 1,2]
and [J--V 1,3].)
Since there is no Frobenius normal form extension for a linear mapping preserving a proper cone
(in a finite-dimensional space),
in this paper we need to adopt a different approach.
In particular,
our treatment of equation (1.1) for the nonnegative matrix case does not rely on a knowledge of the support structure of the generalized eigenvectors of a nonnegative matrix.}
\end{remark}

\begin{remark}
{\rm
Let $A \in \pi (K)$,
and let $\lambda$ be a given positive real number.
By the equivalence of conditions (a) and (b) of Theorem 3.1,
the set $(\lambda I_n -A)K \bigcap K$
(which consists of all vectors $b \in K$ for which equation (1.1) admits a solution)
is equal to $\{ y\in K : \rho_y (A) < \lambda \}$.
The latter set is,
in fact,
an $A$-invariant face of $K$
(see [T--S 2,
the paragraph following Corollary 4.10]).}
\end{remark}

Since a real number $\lambda$ is a distinguished eigenvalue of $A$ for $K$ if and only if $\lambda =\rho_b (A)$ for some $0 \ne b \in K$
(see,
for instance,
[T--W,
Thoerem 2.4(ii)]),
the first part of Remark 3.5 clearly implies the following:

\begin{remark}
{\rm
Let $A \in \pi (K)$,
and let $\lambda$ be a given real number.
Then

(i)
$(\lambda I_n -A)K \bigcap K =K$ if and only if $\lambda > \rho(A)$.

(ii)
$(\lambda I_n -A) K \bigcap K =\{0\}$ if and only if $\lambda$ is less than or equal to the least distinguished eigenvalue of $A$ for $K$.}
\end{remark}

It is clear that Remark 3.6(ii),
in turn,
implies the nontrivial part
(i.e.,
the ``only if" part)
of the known result
([Tam 1,
Theorem 5.1,
(a)$\, \Longleftrightarrow \,$(b)]) that if $A \in \pi (K)$,
then $\rho (A)$ is the only distinguished eigenvalue of $A$ if and only if for any $x \in K$,
$\rho (A)x \ge^K Ax$ implies $\rho(A)x =Ax$.\\

In [F--N 1,
Lemma 1(b)] it is shown that if $A$ is a linear continuous operator on a Banach space $E$,
then for any $\lambda \in \IC\, , ~x, b \in E$,
if $(\lambda I -A)x=b$,
then $\rho_b (A) \leq \rho_x (A) \leq \max \{ |\lambda|, \rho_b (A) \}$.
The proof given in [F--N 1] relies on a use of the concept of local resolvent function.
In our next remark,
we show that in the finite-dimensional case,
we can obtain a slightly stronger conclusion.

\begin{remark}
{\rm
Let $A \in {\cal M}_n (\IC), 0 \neq b\, , ~x \in \IC^n$ and $\lambda \in \IC$ be such that $(\lambda I_n -A)x =b$.
Then we have either $\rho_x (A) > \rho_b(A)$
(in which case $\rho_x (A) =|\lambda|$),
or $\rho _x (A)=\rho_b(A)$
(in which case $\rho_x (A)$ may be greater than,
less than,
or equal to $|\lambda|$.)}
\end{remark}

To show this,
let $x=x_1+\cdots+x_k$ be the representation of $x$ as a sum of generalized eigenvectors of $A$,
and let $\lambda_1,\ldots,\lambda_k$ be the corresponding distinct eigenvalues.
Then the representation of $b$ as a sum of generalized eigenvectors of $A$ consists of the nonzero terms in the sum $(\lambda I_n -A) x_1 +\cdots+(\lambda I_n -A)x_k$.
Now,
for each $j=1, \ldots, k$,
$(\lambda I_n-A)x_j$ is a generalized eigenvector of $A$ corresponding to $\lambda_j$ with the same order as $x_j$,
unless $\lambda=\lambda_j$,
in which case if it is nonzero then its order as a generalized eigenvector is one less than that of $x_j$.
So,
by the definition of local spectral radius,
we have $\rho_b (A)=\rho_x (A)$,
unless there exists some $i$ such that $\lambda=\lambda_i, x_i$ is an eigenvector,
$|\lambda_i|=\rho_x (A)$,
and $|\lambda_j|< \rho_x (A)$ for all $j \neq i$.
In the latter case,
clearly we have $|\lambda|=\rho_x (A) > \rho_b (A)$.\\

By the following remark,
the solution $x^0$ of equation (1.1)
(as described in Theorem 3.1) satisfies $\rho_{x^0} (A) =\rho_b (A)$,
and if $x \in K$ is a solution different from $x^0$ then $\rho_x (A)=\lambda$.

\begin{remark}
{\rm
Let $A \in {\cal M}_n (\IC), 0 \neq b \in \IC^n$ and $\lambda \in \IC$.
If $| \lambda |> \rho_b (A)$,
then $x^0=\sum^{\infty}_{j=0} \lambda^{-j-1} A^j b$ is a solution of the equation $(\lambda I_n -A)x=b$ such that $\rho_{x^0} (A)=\rho_b (A)$.
In this case,
if $x$ is any other solution,
then $\rho_x (A)=| \lambda|$.}
\end{remark}

When $| \lambda| > \rho_b(A)$,
it is clear that $x^0$ is a solution of the equation $(\lambda I_n -A)x=b$.
Since $x^0\in W_b$,
we have $\rho_{x^0} (A) \leq \rho_b (A)$.
By Remark 3.7 we also have $\rho_{x^0} (A) \geq \rho_b (A)$ and hence the equality.
If $x$ is another solution of the said equation,
then we have $(\lambda I -A)(x-x^0)=0$.
In this case,
we can write $x=x^0 +(x-x^0)$,
where $\rho_{x^0} (A) <|\lambda|$ and $x-x^0$ is an eigenvector of $A$ corresponding to $\lambda$.
Hence,
we have $\rho_x (A) =| \lambda|$.

We call a proper cone $K$ {\it subpolar} if $K \subseteq K^*$.
It is clear that every self-dual cone,
in particular,
the nonnegative orthant $\IR^n_+$ is subpolar.
The following is also a corollary of Theorem 3.1.

\begin{coro}
Let $K$ be a subpolar proper cone.
Let $A \in \pi (K)$,
and let $0 \ne b \in K$.
For any positive real number $\lambda$,
the equations
\[
(\lambda I_n -A)x =b , ~ x \in K
\]

and
\[
(A^T -\lambda I_n)z=b, ~ z \in K^*.
\]
cannot be solvable simultaneously.
The same is true for the equations
\[
(A-\lambda I_n)x =b, ~ x \in K
\]

and
\[
(\lambda I_n -A^T)z =b, ~ z \in K^* .
\]
\end{coro}

\noindent
{\it Proof}.
Consider the first two equations.
Assume to the contrary that there exist $x \in K$ and $z \in K^*$ such that
$(\lambda I_n -A) x=b$ and $(A^T -\lambda I_n)z =b$.
Then we have
\[
0 \ge -b^T x =z^T(\lambda I_n -A)x =z^T b \ge 0,
\]
and hence $b^T x =0$.
On the other hand,
by the last part of Theorem 3.1 we can write $x$ as $x^0 +u$,
where $x^0 =\sum^\infty_{j=0} \lambda^{-j-1} A^j b$ and $u$ is either the zero vector or is a distinguished eigenvector of $A$ (for $K$) corresponding to $\lambda$.
But then we have $b^T x \ge b^T x^0 >0$,
which is a contradiction.

Now suppose the last two equations are solvable simultaneously.
By a similar argument we can show that $z^T b$
$(=b^T x)=0$.
By Theorem 3.1 again we can write $z$ as $z^0 +w$,
where $z^0 =\sum^\infty_{j=0} \lambda^{-j-1} (A^T)^j b$ and $w$ is either the zero vector or is a distinguished eigenvector of $A^T$ (for $K^*$);
hence we have $z^T b >0$,
which is a contradiction.\hfill$\Box$\\

\section{The equation {\boldmath $(A-\lambda I_n)x =b\, , ~x \in K$}.}

Given $A \in \pi (K)$ and $\lambda >0$,
it is readily checked that the set $(A-\lambda I_n)K \bigcap K$,
which consists of all vectors $b \in K$ for which equation (1.2) is solvable,
is an $A$-invariant subcone of $K$.
Unlike the cone $(\lambda I_n-A)K \bigcap K$
(for equation (1.1)),
in general,
$(A-\lambda I_n)K \bigcap K$ is not a face of $K$.
(For an example of the latter assertion in the nonnegative matrix case,
consider the nonnegative matrix associated with the singular $M$-matrix given in [T--W,
Example 4.6].)
Moreover,
the question of solvability of equation (1.2) is more delicate.
As we shall see,
it depends on whether $\lambda$ is greater than,
equal to,
or less than $\rho_b(A)$.

\begin{theorem}
Let $A \in \pi (K)$,
let $0 \neq b \in K$,
and let $\lambda$ be a given positive real number such that $\lambda > \rho_b (A)$.
Then the equation {\rm (1.2)} is solvable if and only if $\lambda$ is a distinguished eigenvalue of A for K and $b \in \Phi({\cal N}(\lambda I_n -A) \bigcap K)$.
In this case,
for any solution $x$ of {\rm (1.2)} we have ${\rm sp}_A (x)=(\lambda ,1)$.
\end{theorem}

\noindent
{\it Proof.}
Since $\lambda > \rho_b (A)$,
the vector $\sum^{\infty}_{k=0} \lambda^{-k-1} A^k b$ is well-defined.
We denote it by $x^0$.
Clearly $x^0 \in K$ and $(\lambda I_n -A)x^0=b$.

``If$\,$" part:
Take any vector $u \in {\rm relint}({\cal N}(\lambda I_n -A) \bigcap K)$.
Then $\Phi ({\cal N}(\lambda I_n -A) \bigcap K)$ equals $\Phi (u)$ and is an $A$-invariant face,
as $u$ is a distinguished eigenvector of $A$.
But $b \in \Phi (u)$,
by the definition of $x^0$,
it follows that we also have $x^0 \in \Phi (u)$ and hence $-x^0 \in {\rm span}\,\Phi (u)$.
Choose $\alpha >0$ sufficiently large so that $\alpha u -x^0 \in K$.
Then
\[
(A-\lambda I_n) (\alpha u -x^0)=(A-\lambda I_n)(\alpha u)+(A-\lambda I_n)(-x^0)=b.
\]
Hence $x=\alpha u -x^0$ is the desired solution vector for equation (1.2).

``Only if$\,$" part:
Let $x$ satisfy equation (1.2).
Then we have $(A-\lambda I_n)(x+x^0)=b+(-b)=0$.
Hence $\lambda$ is a distinguished eigenvalue of $A$ and $w=x+x^0$ is a corresponding distinguished eigenvector.
In view of the definition of $x^0$,
we have
\[
b \in \Phi (x^0) \subseteq \Phi(w) \subseteq \Phi ({\cal N}(\lambda I -A ) \bigcap K),
\]
as desired.
Note that by Remark 3.8 we have $\rho_{x^0} (A)=\rho_b(A)$.
But $x=w-x^0$,
$w$ is an eigenvector of $A$ corresponding to $\lambda$,
and $\lambda > \rho _b (A)$,
it follows that we have $\rho_x (A)=\lambda$ and ${\rm ord}_A (x) =1$.\hfill $\Box$\\

Specialized to the nonnegative matrix case,
Theorem 4.1 yields the following:

\begin{coro}
Let P be an $n\! \times \! n$ nonnegative matrix,
let $0 \ne b \in \IR^n_+$,
and let $\lambda$ be a positive real number such that $\lambda > \rho_b (P)$.
Then the equation
\[
(P-\lambda I_n)x=b, ~~x \geq 0
\]
is solvable if and only if $\lambda$ is a distinguished eigenvalue of P such that for any class $\alpha$ of P,
if $\alpha \bigcap {\rm supp}(b)\neq \emptyset$,
then $\alpha$ has access to a distinguished class of P associated with $\lambda$.
\end{coro}

\noindent
{\it Proof.} ``Only if$\,$" part:
By Theorem 4.1 $\lambda$ is a distinguished eigenvalue of $P$ and $b \in \Phi({\cal N}(\lambda I_n -P) \bigcap \IR^n_+ )$.
Let $\alpha_1 , \ldots, \alpha_r$ be the distinguished classes of $P$ associated with $\lambda$.
By the Frobenius-Victory Theorem (see [Sch 3,
Theorem 3.1] or [T--S 1,
Theorem 2.1]),
for each $j=1, \ldots, r$,
there is a (up to multiple)
unique nonnegative eigenvector $x^j=(\xi_1, \ldots, \xi_n)^T$ of $P$ corresponding to $\lambda$ such that $\xi_i >0$ if and only if $i$ has access to $\alpha_j$.
Furthermore,
each vector in the cone ${\cal N}(\lambda I_n -P) \bigcap \IR^n_+$ is a nonnegative linear combination of the vectors $x^1,\ldots,x^r$.
As a consequence,
the $P$-invariant face $\Phi ({\cal N}(\lambda I_n-P)\bigcap \IR^n_+)$ is equal to $\Phi (x^1+\cdots+x^r)$,
and hence can be expressed as $F_I$,
where $I$ is the initial subset for $P$ for which the corresponding initial collection $\cal C$ of classes of $P$ consists of all classes having access to at least one of the ${\alpha_j }'s, ~j=1, \ldots, r$.
Since $b \in \Phi({\cal N}(\lambda I_n -P) \bigcap \IR^n_+)$,
supp($b$) is included in $I$.
If $\alpha$ is a class of $P$ such that $\alpha \bigcap {\rm supp}(b) \neq \emptyset$,
then certainly $\alpha$ has access to a class final in $\cal C$.
But the classes final in $\cal C$ are precisely $\alpha_1, \ldots, \alpha_r$.
Hence our assertion follows.

``If$\,$" part:
Let $I$ still denote the initial subset for $P$ with the same meaning as introduced above.
Then our assumption on supp($b$) clearly implies that supp$(b) \subseteq I$.
Hence,
$b \in F_I =\Phi ({\cal N}(\lambda I_n-P) \bigcap \IR^n_+)$.
By Theorem 4.1 it follows that the equation $(P-\lambda I_n)x=b, ~x \geq 0$ is solvable.\hfill $\Box$\\

Consider equation (1.2) for the case $\lambda > \rho(A)$.
Clearly we have $\lambda > \rho_b(A)$,
so Theorem 4.1 is applicable.
If the equation is solvable,
then necessarily $\lambda$ is a distinguished eigenvalue of $A$ for $K$,
which is a contradiction,
as $\lambda > \rho(A)$.
This proves the following:

\begin{coro}
Let $A\in \pi (K)$,
and let $0 \neq b \in K$.
A necessary condition for equation $(1.2)$ to have a solution is that $\lambda \leq \rho (A)$.
\end{coro}

An alternative way to establish Corollary 4.3 is to make use of the known fact that $\sup \Omega =\rho (A)$ (see [T--W,
Theorem 3.1]).\\

Suppose $(A-\lambda I_n)x=b$,
where $A\in {\cal M}_n (\IC),~0 \neq b, ~x\in \IC^n$ and $\lambda \in \IC$.
Then we have $(\lambda I_n -A)x=-b$,
and by Remark 3.7 we have either $\rho_x(A) > \rho_b (A)$
(in which case $\rho_x (A)=|\lambda|$),
or $\rho_x(A)=\rho_b (A)$
(in which case $\rho_x(A)$ may be greater than,
less than,
or equal to $|\lambda|$).
Our next result implies that when $A \in \pi (K),~0\neq b,~x \in K$ and $\lambda >0$ satisfy $(A-\lambda I_n)x=b$,
the possibility $\rho_x(A)=\rho_b (A) <\lambda$ cannot occur.

\begin{remark}
{\rm
Suppose $(A-\lambda I_n)x=b$,
where $A\in \pi (K),~x,b$ are nonzero vectors of $K$ and $\lambda>0$.
If $\lambda \leq \rho_b (A)$,
then $\rho_x (A)=\rho_b (A)$.
If $\lambda > \rho_b(A)$,
then $\lambda=\rho_x(A)$.}
\end{remark}

It suffices to show that only one of the following three possibilities can occur:
$\lambda=\rho_x (A)=\rho_b(A)$ or $\lambda =\rho_x (A) > \rho_b(A)$ or $\rho_x (A)=\rho_b(A)>\lambda$.
Let $x=x_1+\cdots +x_k$ be the representation of $x$ as a sum of generalized eigenvectors of $A$,
where $\lambda_1,\ldots,\lambda_k$ are the corresponding distinct eigenvalues.
Then the representation of $b$ as a sum of generalized eigenvectors of $A$ consists of the nonzero terms in the sum $(A-\lambda I_n)x_1 +\cdots +(A-\lambda I_n)x_k$,
and by the argument given in the paragraph following Remark 3.7 we see that if $\lambda =\rho_x(A)$,
then $\rho_b(A) \leq \rho_x(A)$,
and if $\lambda \neq \rho_x(A)$,
then $\rho_b (A)=\rho_x(A)$ and ord$_A (b)={\rm ord}_A (x)$.
It remains to show that when $\lambda \neq \rho_x(A)$ we must have $\rho_x(A) >\lambda$.
Since $A \in \pi (K)$ and $0 \ne x \in K$,
by the local Perron-Schaefer condition on $A$ at $x$,
we may assume that $\lambda_1=\rho_x(A)$.
Denote the common value of ord$_A (b)$ and ord$_A(x)$ by $m$.
Applying [T--S 2,
Corollary 4.8]
to the vectors $x$ and $b$ respectively,
we find that the vectors $(A-\rho_x (A) I_n)^{m-1} x_1$ and $(A-\rho_x(A) I_n)^{m-1} (A-\lambda I_n)x_1$ are both distinguished eigenvectors of $A$ corresponding to $\rho_x(A)$.
But then the latter vector is equal $\rho_x (A) -\lambda$ times the former vector.
Hence we have $\rho_x (A)>\lambda$.

\begin{theorem}
Let $A\in \pi (K)$,
and let $0 \ne b \in K$.
If the linear equation
$$(A-\rho_b(A) I_n)x=b\, ,~~x \in K$$
is solvable,
then $b \in (A -\rho_b (A) I_n) \Phi ({\cal N}( (\rho_b(A) I_n -A)^n) \bigcap K)$.
\end{theorem}

\noindent
{\it Proof.}
Suppose that there exists a (nonzero) vector $x\in K$ such that
\setcounter{equation}{0}
\begin{equation}
(A-\rho_b (A) I_n)x=b.
\end{equation}
By Remark 4.4 we have $\rho_x (A)=\rho_b(A)$.
Rewriting (4.1),
we obtain $A x=\rho_b (A) x +b$,
and so $b \in \Phi (A x) \subseteq \Phi (\hat{x})$,
where $\hat{x} =(I_n +A)^{n-1} x$.
By Lemma 2.1 $\Phi (\hat{x})$ is an $A$-invariant face,
and $\rho_{\hat{x}} (A) =\rho_x (A) =\rho_b (A)$.
Equation (4.1) also implies that
\[
(A-\rho_b (A) I_n) \hat{x} = (I_n +A)^{n-1} b \in \Phi(\hat{x}),
\]
i.e.,
$\rho_b (A) \in \Omega_1 (A|_{{\rm span}\Phi (\hat{x})})$.
But $\rho_b (A)=\rho_{\hat{x}} (A)=\rho (A|_{{\rm span}\Phi (\hat{x})} )$,
by [Tam 1,
Theorem 5.1] it follows that $A$ has a generalized eigenvector $y$
(corresponding to $\rho_b (A)$)
that lies in relint $\Phi (\hat{x})$.
Therefore,
we have
$$\Phi (\hat{x}) =\Phi (y) \subseteq \Phi ({\cal N}((\rho_b (A) I_n -A)^n) {\textstyle\bigcap} K).$$
But $b=(A -\rho_b(A) I_n)x$ and $x \in \Phi (\hat{x})$,
so our assertion follows.\hfill$\Box$\\

We can also reformulate our above results in another way:

\begin{theorem}
Let $A \in \pi(K)$,
and let $\lambda$ be a distinguished eigenvalue of $A$ for $K$.
Then
\begin{eqnarray*}
& & (A -\lambda I_n)K {\textstyle \bigcap} \{b \in K: \,\rho_b (A) \le \lambda \} \\
&=& K {\textstyle \bigcap} (A-\lambda I_n) \Phi ({\cal N}((\lambda I_n -A)^n)
      {\textstyle \bigcap} K) \\
&\subseteq & \Phi ({\cal N}((\lambda I_n -A)^n) {\textstyle \bigcap} K) {\textstyle \bigcap}
             \{ y \in K: \, {\rm sp}_A(y) \preceq (\lambda, m_\lambda -1)\},
\end{eqnarray*}
where $m_\lambda$ denotes the maximal order of distinguished generalized eigenvectors of $A$ corresponding to $\lambda$.
\end{theorem}

\noindent
{\it Proof}.
For convenience,
denote the sets $(A -\lambda I_n)K \bigcap \{ b \in K:\, \rho_b(A) \le \lambda\}$,
$K\bigcap (A -\lambda I_n) \Phi ({\cal N}((\lambda I_n -A)^n) \bigcap K)$ and
$\Phi ({\cal N}((\lambda I_n -A)^n)\bigcap K)\bigcap \{y \in K: \, {\rm sp}_A(y)
\preceq (\lambda, m_\lambda -1)\}$ respectively by $S_1$,
$S_2$ and $S_3$.

Let $b\in (A-\lambda I_n)K \bigcap K$.
If $\rho_b(A)=\lambda$,
then by Theorem 4.5,
$b \in S_2$.
If $\rho_b (A) < \lambda$,
then from the ``if" part of the proof of Theorem 4.1,
we also have $b \in (A-\lambda I_n) \Phi ({\cal N}(\lambda I_n -A) \bigcap K) \subseteq S_2$.
This proves $S_1 \subseteq S_2$.

Note that the inclusion $S_2 \subseteq S_1$ follows from $S_2 \subseteq S_3$.
It remains to show the latter inclusion.
Let $b \in S_2$.
Then $b \in K$ and there exists some vector $x \in \Phi({\cal N}((\lambda I_n -A)^n )\bigcap K)$ such that $(A-\lambda I_n)x=b$.
>From the latter equation,
we obtain
\[
b \in \Phi (Ax) \subseteq \Phi (\hat{x}) \subseteq \Phi({\cal N}((\lambda I_n -A)^n){\textstyle\bigcap} K),
\]
where $\hat{x}$ has the same meaning as given in Lemma 2.1,
and the last inclusion holds because $\Phi({\cal N}((\lambda I_n -A)^n)\bigcap K)$ is an $A$-invariant face
(as ${\cal N}((\lambda I_n -A)^n)\bigcap K$ is an $A$-invariant set)
that contains $x$.
If $\rho_b(A) < \lambda$,
we already have $b \in S_3$.
So suppose $\rho_b (A) \ge \lambda$.
Since $b =(A-\lambda I_n)x$,
we have $\rho_b(A) \le \rho_x(A) \le \lambda$,
where the last inequality holds as $x \in \Phi ({\cal N}((\lambda I_n -A)^n)\bigcap K)$
(see [T--S 2,
Theorem 4.9(ii)(b)).
Hence,
we have $\rho_b(A)=\lambda =\rho_x (A)$.
It follows that we also have ord$_A(b)={\rm ord}_A(x)-1$
(see the paragraph following Remark 4.4).
By an argument given in the proof of Theorem 4.5,
we also see that there exists a generalized eigenvector $y$
(corresponding to $\lambda$)
that lies in relint $\Phi(\hat{x})$.
Since $\hat{x}$ and $y$ generate the same face of $K$,
by [T--S 2,
Theorem 4.7] and Lemma 2.1 we have ord$_A(y)={\rm ord}_A(\hat{x})={\rm ord}_A(x)$.
But ord$_A(y) \le m_\lambda$,
so we have ord$_A(b) \le m_\lambda -1$.
This proves that $b \in S_3$.
The proof is complete.\hfill$\Box$\\

We would like to mention that the last of the three sets considered in the preceding theorem is an $A$-invariant face
(see [T--S 2,
Corollary 4.10]),
and also in general the face of $K$ generated by the first
(and hence also the second) set can be strictly included in the last set.
One can readily find an illustrative example in the nonnegative matrix case.
As for the quantity $m_\lambda$,
we would like to point out the following:

\begin{remark}
{\rm
Let $\lambda$ be a distinguished eigenvalue of $A$ for $K$.
Let $m_{\lambda}$ denote the maximal order of distinguished generalized eigenvectors of $A$ corresponding to $\lambda$.
Also let $G_{\lambda}$ denote the $A$-invariant face $\{ y \in K:\, \rho_y(A) \le \lambda\}$ of $K$.
In general,
we have $m_{\lambda} \le \nu_{\lambda}(A|_{G_{\lambda}}) \le \nu_{\lambda}(A)$,
where both inequalities can be strict.
The first inequality becomes an equality when the cone $K$ is polyhedral,
whereas the second inequality becomes an equality when $\lambda=\rho(A)$.}
\end{remark}

\begin{coro}
For any $A\in \pi (K)$,
we have
\[
\begin{array}{cl}
 & \Phi(( A-\rho(A) I_n)K \bigcap K) \\
=& \Phi(K \bigcap [(A-\rho(A) I_n) \Phi({\cal N}((\rho(A) I_n-A)^n) \bigcap K)]) \\
\subseteq & \Phi({\cal N}((\rho(A)I_n -A)^n)\bigcap K)\bigcap d_{K^*}(\Phi({\cal N}(\rho(A)I_n-A^T) \bigcap K^*)). \\
\subseteq & \Phi({\cal N}((\rho(A)I_n-A)^n)\bigcap K)\bigcap \{ y \in K: \, {\rm sp}_A (y) \preceq (\rho(A), \nu_{\rho(A)} (A)-1) \}.
\end{array}
\]
\end{coro}

\noindent
{\it Proof}.
Applying Theorem 4.6 to the case $\lambda =\rho(A)$,
we readily see that the first two sets are equal and are included in $\Phi({\cal N}((\rho(A)I_n -A)^n)\bigcap K)$.
Now consider any $b \in (A-\rho(A) I_n)K \bigcap K$.
Then $b =(A-\rho(A) I_n)x$ for some $x \in K$.
Take any $z \in {\cal N}(\rho(A) I_n -A^T) \bigcap K^*$.
We have
\[
\langle b, \, z \rangle =\langle (A-\rho(A)I_n)x , \, z \rangle = \langle x, \, (A^T-\rho(A)I_n)x \rangle =0,
\]
hence
$\langle b^\prime,~z^\prime \rangle = 0$ for any choice of
$b^\prime \in \Phi ((A-\rho(A) I_n)K \bigcap K)$ and
$z^\prime \in \Phi({\cal N}(\rho(A) I_n -A^T) \bigcap K^*)$.
This shows that $\Phi((A-\rho(A)I_n)K \bigcap K)\subseteq d_{K^*}$
$(\Phi({\cal N}(\rho(A)I_n -A^T)\bigcap K^*)$,
and hence the first inclusion.

To establish the second inclusion,
it suffices to prove $d_{K^*} (\Phi ({\cal N}(\rho(A)I_n -A^T) \bigcap K^*) \subseteq \{ y \in K: {\rm sp}_A(y) \preceq (\rho(A), \nu -1)\}$,
where for convenience we use $\nu$ to denote $\nu_{\rho(A)} (A)$.
Consider any vector $x\in d_{K^*}(\Phi({\cal N}(\rho(A) I_n -A^T) \bigcap K^*)$.
Denote by $E^{(0)}_\rho (A)$ the projection of $\IC^n$ onto ${\cal N}((\rho(A) I_n -A)^n)$ along ${\cal R}((\rho(A) I_n -A)^n)$.
Then,
according to [Sch 2,
Theorem 5.2],
the restriction of $(A-\rho(A)I_n)^{\nu -1} E^{(0)}_\rho (A)$ to $\IR^n$ is nonzero and belongs to $\pi (K)$.
A similar assertion can be said of $(A^T-\rho(A)I_n)^{\nu -1} E^{(0)}_\rho (A^T)$,
where $E^{(0)}_\rho (A^T)$ has a similar meaning.
Now take any $z \in \mbox{int } K^*$.
Then $(A^T -\rho(A) I_n)^{\nu -1} E^{(0)}_\rho (A^T)z \in K^*$ and is an eigenvector of $A^T$ corresponding to $\rho(A)$.
By our choice of $x$,
we have
\[
\begin{array}{ccl}
0 &=& \langle x, \, (A^T-\rho(A) I_n)^{\nu -1} E^{(0)}_\rho (A^T) z \rangle \\
  &=& \langle (E^{(0)}_\rho (A^T))^T (A-\rho(A) I_n)^{\nu -1} x, \, z \rangle \\
  &=& \langle (A-\rho(A) I_{n-1} )^{\nu -1} E^{(0)}_\rho (A) x, \, z \rangle ,
\end{array}
\]
where the last equality follows from the fact that $(E^{(0)}_\rho (A^T))^T=E^{(0)}_\rho (A)$,
and $(A-\rho(A)I_{n-1})^{\nu -1}$ and $E^{(0)}_\rho (A)$ commute.
But $z \in \mbox{int } K^*$,
so we have $(A-\rho(A) I_{n-1})^{\nu -1} E^{(0)}_\rho (A) x =0$,
which is equivalent to that sp$_A(x) \preceq (\rho(A), \, \nu-1)$.
The proof is complete.\hfill$\Box$\\

By applying Corollary 4.8 to $A|_{\mbox{span }F} \in \pi (F)$,
where $F=\{b \in K: \rho_b (A) \le \lambda \}$,
one can readily obtain an $A$-invariant face that includes the second set but is included in the third set of Theorem 4.6.
We leave it to the reader to write out the set.

In words,
the third set that appears in Corollary 4.8 is the intersection of the following two $A$-invariant faces of $K$:
the face generated by the distinguished generalized eigenvectors of $A$ corresponding to $\rho(A)$,
and the dual face of the face of $K^*$ generated by the distinguished eigenvectors of $A^T$.
As the following example will show,
in general,
the third set in Corollary 4.8 can strictly include the second
(and hence also the first) set.
The example is borrowed from [T--W,
Example 3.7].

\begin{example}
{\rm
Let $\alpha$ be a given real number with $0 < \alpha < 1$.
Let $C$ be the closed convex set in $\IR^2$ with extreme points $(k, \alpha^{k-1})^T$,
$k=1, 2, \ldots$ and with recession cone $O^+ C=\{\lambda (1,0)^T: \lambda \ge 0\}$.
Let $K$ be the cone in $\IR^3$ given by
\[
K= \left\{ \lambda {x \choose 1} \in \IR^3: x \in C, \, \lambda \ge 0 \right\} \bigcup
\{ \lambda (1, 0, 0)^T: \lambda \ge 0 \}.
\]
Also let
\[
A= \left[ \begin{array}{ccc}
1 & 0 & 1 \\
0 & \alpha & 0 \\
0 & 1 & 1
\end{array} \right] .
\]
Then $K$ is a proper cone in $\IR^3$ and we have $A \in \pi (K)$.
Clearly,
the face $\Phi({\cal N}((\rho(A) I_3 -A)^3 )\bigcap K)$ is equal to the extreme ray $\Phi ((1, 0, 0)^T)$ of $K$.
By Corollary 4.8 it follows that $\Phi ((A-\rho(A) I_3)K \bigcap K)=\{0\}$.
(By solving the equation $(A-\rho(A) I_3)x =b,~x \in K$,
where $b \in K$,
one can also readily show that the set $(A-\rho(A) I_3) K \bigcap K$ equals $\{0\}$.)
On the other hand,
it is not difficult to show that $d_{K^*} (\Phi ({\cal N}(\rho(A) I_3 -A^T)\bigcap K^*))=\Phi((1,0,0)^T)$.
Hence,
the third set of Corollary 4.8 equals $\Phi ((1, 0,0)^T)$ and strictly includes the first
(and hence also the second) set.}
\end{example}

Corollary 4.8 also has the following interesting consequence.

\begin{coro}
Let $A \in \pi (K)$.
If $A$ has no distinguished generalized eigenvectors,
other than distinguished eigenvectors,
corresponding to $\rho(A)$,
and if $A$ has no eigenvectors in $\Phi ({\cal N}((\rho(A)I_n -A)^n)\bigcap K)$ corresponding to an eigenvalue other than $\rho(A)$,
then $(A-\rho(A) I_n)K \bigcap K=\{0\}$.
The converse also holds if the cone $K$ is polyhedral.
\end{coro}

\noindent
{\it Proof}.
For convenience,
we denote the set ${\cal N}((\rho(A) I_n -A)^n) \bigcap K$ by $S$.
By our assumption,
it is clear that $S$ is equal to ${\cal N}(\rho(A)I_n -A) \bigcap K$.
Let $b \in K \bigcap (A-\rho(A)I_n)\Phi (S)$.
Then there exists $x \in \Phi(S)$ such that $(A-\rho(A) I_n)x =b$.
Since $\Phi(S)$ is an $A$-invariant face,
we also have $b \in K \bigcap \mbox{span } \Phi(S)=\Phi(S)$.
By our assumption,
$S$ is equal to ${\cal N}(\rho(A) I_n -A)\bigcap K$.
Hence,
we have $\nu_{\rho(A)} (A|_{\mbox{span }\Phi(S)})=1$ and so $\rho_b (A) < \rho(A)$.
Moreover,
our assumption also implies that $\rho(A)$ is the only distinguished eigenvalue of $A|_{\mbox{span }\Phi(S)}$ for $\Phi(S)$.
So,
necessarily,
we have $b=0$.
This proves that $K \bigcap (A-\rho(A) I_n)\Phi (S)=\{0\}$,
and by Corollary 4.8 it follows that $\Phi((A-\rho(A)I_n)K \bigcap K)$,
and hence $(A-\rho(A) I_n)K \bigcap K$ equals $\{0\}$.
[One can also arrive at $\Phi(S) \bigcap (A-\rho(A) I_n)\Phi(S)=\{0\}$ by applying [Tam 1,
Theorem 5.4,
(a)$\Longleftrightarrow$(b)] to $A|_{\mbox{span }\Phi(S)}$.]

Conversely,
suppose $(A-\rho(A) I_n)K \bigcap K=\{0\}$.
Assume,
in addition,
that $K$ is polyhedral.
Let $S$ have the same meaning as before.
First,
suppose $A$ has a distinguished generalized eigenvector other than eigenvector corresponding to $\rho(A)$,
or equivalently,
$\nu_{\rho(A)} (A|_{\mbox{span }\Phi(S)})\ge 2$.
Since $\Phi (S)$ is polyhedral,
by [Tam 1,
Theorem 7.5] there exists a vector $y \in \Phi(S)$ such that $0 \ne (A-\rho(A) I_n)^i y \in \Phi (S)$ for $i =1, \ldots, \nu_{\rho(A)} (A|_{\mbox{span }\Phi(S)}) -1$.
Then the set $(A-\rho(A)I_n) K \bigcap K$ contains the nonzero vector $(A-\rho(A) I_n)y$,
which is a contradiction.
So we must have $S={\cal N}(\rho(A) I_n -A) \bigcap K$.

Now assume that $A$ has a distinguished eigenvector in $\Phi(S)$,
say $w$,
corresponding to an eigenvalue $\lambda$,
which is different from $\rho(A)$.
Choose any $u \in \mbox{relint }S$.
Then there exists $\alpha >0$ such that $u -\alpha w \in \Phi (S)$,
and we have $(A-\rho(A) I_n) (u - \alpha w) =\alpha (\rho(A) -\lambda)w$,
which is a nonzero vector in $(A-\rho(A) I_n) K \bigcap K$.
So we arrive at a contradiction.
The proof is complete.\hfill$\Box$\\

Later,
in Section 5,
the result of Corollary 4.10 will be superseded by Theorem 5.11.

If the polyhedrality assumption on $K$ is dropped,
the converse part of Corollary 4.10 no longer holds.
Here is a counter-example,
which we borrows from the proof of [T--S 1,
Theorem 7.13].
(For an easier counter-example,
one may also use $A$ and $K$ of [Tam 1, Example 5.5].)

\begin{example}
{\rm
Let $n \ge 3$ be an odd integer,
and let $\lambda$ be a positive real number.
Denote by $C$ the unbounded convex set in $\IR^{n-1}$ with extreme points
\[
x_k =\left( {k \choose n-1}, {k \choose n-2}, \ldots, {k \choose 1}\right)^T,
~k=0, \pm 1, \pm 2, \ldots,
\]
and recession cone $O^+ C={\rm ray}((1,0, \ldots, 0)^T)$.
Let $K$ be the proper cone in $\IR^n$ given by:
\[
K = \left\{ \alpha {x \choose 1}: \alpha \ge 0, ~x \in C \right\} \bigcup (O^+ C \times \{0\}).
\]
Let $J_n(\lambda)$ denote the $n \times n$
(upper triangular) elementary Jordan block corresponding to $\lambda$.
Then it is not difficult to verify that $J_n (\lambda)\in \pi (K)$.
Note that all vectors in $K$,
except for those that lie on the extreme ray $\Phi (e_1)$
(where $e_1$ denotes the vector in $\IR^n$ with 1 at its first component and 0 elsewhere),
have positive last component.
As a consequence,
we have $(J_n(\lambda)-\lambda I_n) K \bigcap K =\{0\}$.
However,
all vectors in $K$,
except for those that lie on $\Phi(e_1)$,
are distinguished generalized eigenvectors of $J_n(\lambda)$,
other than eigenvectors,
corresponding to $\rho(J_n (\lambda))$ $(=\lambda)$. }
\end{example}

The question of when $(A -\rho(A) I_n)K \bigcap K=\{0\}$ is related to the problem of determining all real numbers $\lambda$ for which $(A-\lambda I_n) K \bigcap K =\{0\}$.
Indeed,
in view of the following,
the set of all such $\lambda$ consists of all real numbers greater than $\rho(A)$,
together with or without $\rho(A)$.

\begin{remark}
{\rm
Let $A\in \pi(K)$.
For any real number $\lambda$,
we have
\[
\begin{array}{rcl}
\lambda > \rho(A) &\Longrightarrow & (A-\lambda I_n)K \bigcap K=\{0\} \\
&\Longrightarrow & \lambda \ge \rho(A) .
\end{array}
\]
}
\end{remark}

The preceding remark follows simply from the definition of $\Omega (A)$ and the known fact that sup $\Omega (A)=\rho(A)$.

We have shown that,
in general,
the first inclusion in Corollary 4.8 can be strict.
However,
in the nonnegative matrix case,
we always have an equality.

\begin{theorem}
Let $P$ be an $n \times n$ nonnegative matrix.
Let $I$ denote the union of all classes $\alpha$ of $P$ such that $\alpha >\!\!-\,\beta$ for some basic class $\beta$ of $P$.
Then $\Phi((P-\rho(P)I_n)\IR^n_+ \bigcap \IR^n_+)$ is equal to the $P$-invariant face $F_I$ of $\IR^n_+$.
\end{theorem}

\noindent
{\it Proof}.
Let $\alpha_1, \ldots, \alpha_q$ denote the basic classes (of $P$).
By the Nonnegative Basis theorem
(see [Sch 3,
Theorem 7.1] or [Tam 1,
Theorem 5.2]),
the generalized eigenspace ${\cal N}((\rho(P) I_n -P)^n)$ has a basis consisting of nonnegative vectors $x^1, \ldots, x^q$ such that for $j=1, \ldots, q$,
the $i$th component of $x^j$ is positive if and only if $i$ has access to $\alpha_j$.
Clearly,
$x^1+\cdots+x^q \in {\rm relint}({\cal N}((\rho(P) I_n -P)^n)\bigcap \IR^n_+)$ and supp$(x^1+\cdots+x^q)=J$,
where we use $J$ to denote the initial subset for $P$ which is equal to the union of all classes having access to one of the $\alpha_j$'s,
$j=1, \ldots, q$.
Hence,
we have $\Phi({\cal N}((\rho(P) I_n -P)^n)\bigcap \IR^n_+)=F_J$.
Similarly,
by applying the Frobenius-Victory theorem
(see the proof of the ``only if" part of Corollary 4.2),
we also have $\Phi({\cal N}(\rho(P) I_n -P^T)\bigcap \IR^n_+)=F_L$,
where $L$ denotes the union of all classes $\beta$ such that $\beta$
has access from some basic class
(in the digraph of $P$) which is at the same time a distinguished basic class of $P^T$.
Then $d_{\IR^n_+}(\Phi ({\cal N}(\rho(P)I_n -P^T) \bigcap \IR^n_+))$ is equal to $F_M$,
where $M =\langle n \rangle \backslash L$ is an initial subset for $P$.
Denote by ${\cal C}_I$
(respectively,
${\cal C}_J$ and ${\cal C}_M$)
the initial collection of classes associated with the initial set $I$
(respectively,
$J$ and $M$).
Note that the classes that are final in ${\cal C}_J$ are precisely the basic classes which are at the same time distinguished basic classes of $P^T$.
Using this,
it is not difficult to see that ${\cal C}_I={\cal C}_J \bigcap {\cal C}_M$,
hence $I=J \bigcap M$.
Now by Corollary 4.8 and Theorem 2.2,
we have
\[ \begin{array}{cl}
& \Phi((P-\rho(P)I_n) \IR^n_+ \bigcap \IR^n_+) \\
\subseteq & \Phi({\cal N}((\rho(P)I_n -P)^n)\bigcap \IR^n_+) \bigcap
            d_{\IR^n_+}({\cal N}(\rho(P)I_n -P^T) \bigcap \IR^n_+) \\
= & F_J \bigcap F_M \\
= & F_{J \bigcap M} \\
= & F_I .
\end{array}
\]

To prove the reverse inclusion,
we contend that for any distinguished basic class $\alpha$ of $P^T$,
there exists a vector $b^\alpha \in \IR^n_+$ such that the subvector $(b^\alpha)_\beta$ is nonzero if $\beta >\!\!- \, \alpha$ and is zero,
otherwise,
for which the equation $(P-\rho(P)I_n)x=b^\alpha$,
$x \in \IR^n_+$ is solvable.
Once the contention is proved,
we set $b=\sum b^\alpha$,
where the summation runs through all distinguished basic classes $\alpha$ of $P^T$.
It is not difficult to see that the smallest initial subset for $P$ that includes
supp$(b)$ is $I$.
By Theorem 2.2 this means that $F_I$ is the smallest $P$-invariant face that contains $b$.
But,
by our choice of $b$,
clearly the equation $(P-\rho(P)I_n)x=b$,
$x \in K$ is solvable;
that is,
$b$ also belongs to the $P$-invariant face $\Phi((P-\rho(P)I_n)\IR^n_+ \bigcap \IR^n_+)$.
Thus the reverse inclusion $F_I \subseteq \Phi ((P-\rho(P)I_n)\IR^n_+ \bigcap \IR^n_+)$ also holds.

It remains to prove our contention.
Consider any distinguished basic class $\alpha$ of $P^T$.
We are going to find a pair of vectors $x , ~b \in \IR^n_+$ such that $b_\beta$ is nonzero,
nonnegative if $\beta >\!\!- \, \alpha$ and is zero,
otherwise,
and moreover $b=(P-\rho(P)I_n)x$.
(Here,
for simplicity,
we write the above-mentioned vector $b^\alpha$ as $b$.)
Note that the latter equation amounts to
\begin{equation}
b_\beta =(P_{\beta \beta} -\rho(P)I_n)x_\beta +\sum_{\beta >\!- \, \gamma} P_{\beta \gamma} x_{\gamma}
\end{equation}
for all classes $\beta$ (of $P$).
For any class $\beta$ that has no access to $\alpha$,
we set $x_\beta=b_\beta =0$.
It is easy to check that (4.2) is satisfied for all such classes $\beta$.
Take $x_\alpha$ to be a Perron vector of the irreducible nonnegative matrix
$P_{\alpha \alpha}$ and also set $b_\alpha =0$.
Then (4.2) also holds for $\beta =\alpha$.
Now we use the trace down method to determine $x_\beta$ and $b_\beta$ for
the remaining classes $\beta$.
At a general step,
for a given class $\beta >\!\!- \, \alpha$,
suppose we have already determined the subvectors $x_\gamma$ and $b_\gamma$
for all classes $\gamma$,
$\beta >\!\!- \gamma$,
in such a way that $x_\gamma$ is positive if $\gamma >= \alpha$,
and $b_\gamma$ is nonzero,
nonnegative if $\gamma >\!\!- \, \alpha$,
and $x_\gamma$ (also $b_\gamma$) is zero,
otherwise.
If $\beta$ is a basic class,
we take $x_\beta$ to be a Perron vector of $P_{\beta \beta}$ and $b_{\beta}$ to be the nonzero,
nonnegative vector $\sum_{\beta >\!-\, \gamma} P_{\beta \gamma} x_\gamma$
(noting that there is at least one $\gamma$,
$\beta >\!\!-\, \gamma$,
such that $P_{\beta \gamma} \ne 0$ and $x_\gamma$ is positive).
If $\beta$ is nonbasic,
we take $b_\beta$ to be $\frac{1}{2} \sum_{\beta >\!- \, \gamma} P_{\beta \gamma} x_\gamma$ and $x_\beta$ to be $(\rho(P)I-P_{\beta \beta})^{-1} (\frac{1}{2} \sum_{\beta >\!-\, \gamma} P_{\beta \gamma} x_\gamma)$
(noting that $(\rho(P)I-P_{\beta \beta})^{-1}$ is a positive matrix,
as $\rho(P) > \rho(P_{\beta \beta})$).
In any case,
we have $x_\beta$ is positive,
$b_\beta$ is nonzero,
nonnegative,
and the equation (4.2) is satisfied.
Proceeding in this way,
after a finite number of steps,
we can construct the desired vector $b$ (and $x$).
This establishes our contention.
The proof is complete.\hfill$\Box$\\

We {\it do not know} whether the first inclusion in Corollary 4.8 becomes an equality when the underlying cone $K$ is polyhedral.

Note that Theorem 4.13 tells,
in particular,
that if $P$ is an $n \times n$ nonnegative matrix and if $0 \ne b \in \IR^n_+$ is such that the equation $(P-\rho(P)I_n)x=b$,
$x \ge 0$ is solvable,
then for any class $\alpha$ of $P$ for which $\alpha \bigcap {\rm supp}(b) \ne \emptyset$,
we have $\alpha >\!\!-\, \beta$ for some basic class $\beta$ of $P$.
This strengthens the result of [T--W,
Lemma 4.5],
which is formulated in terms of a singular $M$-matrix,
and in which we have $\alpha >= \beta$ instead of $\alpha >\!\!-\, \beta$ in the conclusion.

For the case of a distinguished eigenvalue,
we have the following:

\begin{theorem}
Let $P$ be an $n \times n$ nonnegative matrix and let $\lambda$ be a distinguished eigenvalue of $P$.
Let $I$ denote the union of all classes $\alpha$ of $P$ such that $\alpha >\!\!-\, \beta$ for some semi-distinguished class $\beta$ of $P$ associated with $\lambda$.
Then $\Phi((P-\lambda I_n)\IR^n_+ \bigcap \{ b \in \IR^n_+: \, \rho_b(P) \le \lambda\} )$ is equal to $F_I$.
\end{theorem}

\noindent
{\it Proof}.
Apply Theorem 4.13 to the principal submatrix $P_{JJ}$ of $P$,
where $J$ is the union of all classes $P$ that have access to some semi-distinguished class of $P$ associated with $\lambda$.\hfill$\Box$\\

As can be readily seen,
Theorem 4.14 implies that for an $n \times n$ nonnegative matrix $P$ and any $0\ne b\in\IR^n_+$,
if the equation $(P-\rho_b(P)I_n)x=b$,
$x \ge 0$ is solvable,
then for any class $\alpha$ of $P$ for which $\alpha \bigcap {\rm supp}(b)\ne \emptyset$,
we have $\alpha >\!\!-\, \beta$ for some semi-distinguished class $\beta$ associated with $\rho_b(P)$.

Now we consider equation (1.2) for $\lambda < \rho_b (A)$.
As we shall see,
in this case,
for the equation to be solvable it is not necessary that $\lambda$ be an eigenvalue of $A$.

We are going to treat the case when $A$ is $K$-irreducible first.
Recall that a matrix $A \in \pi (K)$ is said to be $K$-{\it irreducible} if $A$ leaves invariant no faces of $K$ other than $\{0\}$ and $K$ itself.
Also,
we call a matrix $A$ $K$-{\it positive} if $A (K \backslash \{0\}) \subseteq \mbox{int }K$.

For a square matrix $A$,
we denote by adj $A$ its (classical) adjoint
(i.e.,
the transpose matrix of cofactors).

The following result is known in the standard case of nonnegative matrices.
It first appeared explicitly in [Fan,
Theorem 3] in the context of an $N$-matrix.
(We call an $n \times n$ matrix $A$ and {\it N-matrix} if $A$ is of the form $\lambda I_n -P$,
where $P\ge 0$,
and $\rho_{n-1} (P) < \lambda < \rho(P)$,
where $\rho_{n-1}(P)$ denotes the maximum of the spectral radii of the $(n-1) \times (n-1)$
principal submatrices of $P$.
Note that,
in this case,
the matrix $P$ is necessarily irreducible.)
The proof is almost the same as in the nonnegative matrix case.
(In fact,
parts of the argument can be found in the original papers of Perron [Per 1,2]
and Frobenius [Fro,
1--3].)

\begin{theorem}
If $A\in \pi (K)$ is $K$-irreducible,
then for $\lambda <\rho(A)$,
sufficiently close to $\rho(A)$,
$(A-\lambda I_n)^{-1}$ and {\rm adj}$(\lambda I_n-A)$ are both $K$-positive.
\end{theorem}

\noindent
{\it Proof}.
Let $A \in \pi (K)$ be $K$-irreducible.
Then $\rho(A)$ is a simple eigenvalue and by a standard argument
(see,
for instance,
[B--P,
Corollary 2.2.13]) adj$(\rho(A) I_n-A)$ is a rank-one matrix of the form $yz^T$,
where $y$ and $z$ are eigenvectors of $A$ and $A^T$ respectively,
both corresponding to $\rho(A)$.
Here $y$ or $-y$
(respectively,
$z$ or $-z$) belongs to int $K$
(respectively,
int $K^*$).
Hence,
we have,
either adj$(\rho(A) I_n -A)$ or its negative is $K$-positive.
But
\[
z^Ty=\mbox{tr } (yz^T)=\mbox{tr adj}(\rho(A) I_n -A) = \triangle^\prime (\rho(A)) >0,
\]
where we use $\triangle (t)$ to denote the polynomial det$(tI_n-A)$,
and the last equality follows from the definition of $\rho(A)$ and the fact that $\rho(A)$ is a simple eigenvalue.
Hence,
adj$(\rho(A)I_n -A)$ must be $K$-positive.

[Alternative argument:
Note that if $\lambda > \rho(A)$,
then det$(\lambda I_n -A) >0$ because the polynomial det$(t I_n -A)$ tends to infinity with $t$ and there are no roots larger than $\rho(A)$.
Moreover,
for such $\lambda$,
$(\lambda I_n -A)^{-1}$ is $K$-positive.
Hence,
adj$(\lambda I_n -A)$ is $K$-positive for $\lambda > \rho(A)$.
Letting $\lambda \to \rho(A)$,
we obtain adj$(\rho(A) I_n -A) \in \pi (K)$.
Together with the above,
we can now conclude that adj$(\rho(A) I_n -A)$ is $K$-positive.]

It is not difficult to show that a matrix which is close to a $K$-positive matrix is still $K$-positive.
Hence,
adj$(\lambda I_n -A)$ is $K$-positive for $\lambda < \rho(A)$,
sufficiently close to $\rho(A)$.
Since $\rho(A)$ is a simple eigenvalue,
for $\lambda < \rho(A)$,
sufficiently close to $\rho(A)$,
we have det $(\lambda I_n -A)<0$.
So,
for such $\lambda$,
$(A-\lambda I_n)^{-1}$ is also $K$-positive.\hfill$\Box$\\

The following is an immediate consequence of Theorem 4.15:

\begin{coro}
If $A$ is $K$-irreducible,
then for $\lambda < \rho(A)$,
sufficiently close to $\rho(A)$,
equation $(1.2)$ is solvable for all $b \in K$.
\end{coro}

It is not difficult to show the following:

\begin{remark}
{\rm
Let $A\in \pi (K)$.
For any real number $\lambda$,
the condition that $A-\lambda I_n$ is nonsingular and $(A-\lambda I_n)^{-1} \in\pi (K)$ is equivalent to $(A-\lambda I_n)K \supseteq K$,
or equivalently,
$(A-\lambda I_n)K \bigcap K=K$.
When the equivalent conditions are satisfied,
$\lambda$ must be less than the least distinguished eigenvalue of $A$ for $K$,
as well as the least distinguished eigenvalue of $A^T$ for $K^*$.}
\end{remark}

Now we know that if $A$ is $K$-irreducible,
then for $\lambda < \rho(A)$,
sufficiently close to $\rho(A)$,
we have $K \subseteq (A-\lambda I_n)K$.
On the other hand,
if $\lambda < \rho(A)$ but is not close to $\rho(A)$,
then this is no longer true.
To see this,
consider an $n \times n$ positive matrix $P$ and any real number $\lambda$ less than the minimum of the diagonal entries of $P$.
Then $P-\lambda I_n$ is a positive matrix,
so $(P-\lambda I_n) (\IR^n_+ \backslash \{0\}) \subseteq \mbox{int } \IR^n_+$ and hence
$\IR^n_+ \not\subseteq (P-\lambda I_n)\IR^n_+$.

The above example also shows that the converse of the last part of Remark 4.17 is not true.
Note,
however,
that there are examples of $A \in \pi (K)$ for which it is true that $(A-\lambda I_n)^{-1} \in \pi (K)$ for all real numbers $\lambda$ less than the least distinguished eigenvalue of $A$ for $K$.
For instance,
take $K=\IR^2_+$ and $A=\mbox{diag}(1,2)$.

For $\lambda < \rho(A)$ and a $K$-irreducible matrix $A$,
in general,
the truth is that,
we have $\Phi((A-\lambda I_n)K \bigcap K)=K$.
This follows simply form the fact that in this case the set $(A-\lambda I_n)K \bigcap K$ must contain the Perron vector of $A$,
which necessarily lies in int $K$.

Corollary 4.16 also yields the following result,
true for $A$ that need not be $K$-irreducible.

\begin{coro}
Let $A \in \pi (K)$,
and let $0 \ne b \in K$.
If $\Phi((I_n +A)^{n-1} b)$ is a minimal nonzero $A$-invarinat face of $K$,
then for $\lambda < \rho_b(A)$,
sufficiently close to $\rho_b(A)$,
equation $(1.2)$ has a solution.
\end{coro}

\noindent
{\it Proof}.
For convenience,
we write $(I_n +A)^{n-1} b$ as $\hat{b}$.
By Lemma 2.1,
$\Phi(\hat{b})$ is always an $A$-invariant face.
Apply Corollary 4.16 to $A|_{\mbox{span }\Phi(\hat{b})} \in \pi (\Phi(\hat{b}))$,
noting that $A|_{\mbox{span }\Phi(\hat{b})}$ is irreducible with respect to $\Phi(\hat{b})$,
as $\Phi(\hat{b})$ is minimal nonzero $A$-invariant face of $K$.\hfill$\Box$\\

More generally,
we have the following:

\begin{theorem}
Let $A\in \pi (K)$.
Let $r$ denote the largest real eigenvalue of $A$ less than $\rho(A)$.
$($If no such eigenvalues exist,
take $r=- \infty$.$)$
Then for any $\lambda$,
$r < \lambda < \rho(A)$,
we have
\[
\Phi((A-\lambda I_n)K {\textstyle \bigcap} K)=\Phi ({\cal N}((\rho(A)I_n-A)^n) {\textstyle \bigcap}K).
\]
\end{theorem}

\noindent
{\it Proof}.
For convenience,
denote the $A$-invariant faces $\Phi((A-\lambda I_n)K \bigcap K)$ and $\Phi ({\cal N}((\rho(A)I_n-A)^n) \bigcap K)$ of $K$ respectively by $C_1$ and $C_2$.
It is easy to see that if $y \in K$ is an eigenvector of $A$ corresponding to $\rho(A)$ then $y \in (A-\lambda I_n)K \bigcap K$;
hence we have,
$\rho(A) \ge \rho(A|_{{\rm span}\, C_1}) \ge \rho_y (A)=\rho(A)$,
i.e.,
$\rho(A|_{{\rm span}\, C_1})=\rho(A)$.
For $i=1,2$,
denote by $C^D_i$ the dual of $C_i$ in span$\,C_i$ and by
Ad$(A|_{{\rm span}\, C_i})$ the adjoint of the linear operator $A|_{{\rm span}\, C_i}$.
Suppose Ad$(A|_{{\rm span}\, C_1})$ has a distinguished eigenvalue for $C_1^D$
other than $\rho(A)$,
say $\mu$,
and let $z\in C^D_1$ be the corresponding eigenvector.
By our choice of $\lambda$,
it is clear that $\mu < \lambda$.
Choose any vector $b \in {\rm relint}((A-\lambda I_n)K \bigcap K)$,
and let $x \in K$ be a vector such that $(A-\lambda I_n)x=b$.
Then,
in fact,
we have $x \in {\rm span} \, C_1 \bigcap K =C_1$
(see the discussion following Remark 4.13).
Also,
$\langle z, b \rangle > 0$ as $z \in C_1^D$ and $b \in {\rm relint}\, C_1$.
On the other hand,
we have
\begin{eqnarray*}
\langle z, b \rangle &=& \langle z, (A|_{{\rm span}\, C_1} -\lambda I)x \rangle \\
&=& \langle {\rm Ad}(A|_{{\rm span}\, C_1} -\lambda I) z, x \rangle \\
&=& (\mu -\lambda ) \langle z, x \rangle \\
&\le& 0 ,
\end{eqnarray*}
as $\mu < \lambda , z \in C_1^D$ and $x \in C_1$.
So we arrive at a contradiction.
This proves that $\rho(A)$ is the only distinguished eigenvalue of Ad$(A|_{{\rm span}\, C_1})$ (for $C_1^D$).
By [Tam 1,
Theorem 5.1] this implies that relint$\, C_1$ contains a generalized eigenvector of $A$ corresponding to $\rho(A)$.
It follows that we have $C_1 \subseteq C_2$.

To prove the reverse inclusion,
consider $A|_{{\rm span}\, C_2} \in \pi (C_2)$.
Clearly,
relint$\, C_2$ contains a generalized eigenvector of $A|_{{\rm span}\, C_2}$.
By [Tam 1,
Theorem 5.1] again,
$\rho(A|_{{\rm span}\, C_2})$,
which is $\rho(A)$,
is the only distinguished eigenvalue of Ad$(A|_{{\rm span}\, C_2})$ for $C_2^D$.
But according to [T--W,
Theorem 3.3],
$\sup$ $\Omega_1 (A|_{{\rm span}\, C_2})$ is equal to the least distinguished eigenvalue of Ad$(A|_{{\rm span}\, C_2})$;
so we have $\sup \, \Omega_1(A|_{{\rm span}\, C_2})=\rho(A)$.
Since $\lambda < \rho(A)$,
we can find a $\lambda^\prime$,
$\lambda < \lambda^\prime \le \rho(A)$,
such that $\lambda^\prime \in \Omega_1(A|_{{\rm span}\, C_2})$.
Then there exists $u \in {\rm relint}\, C_2$ such that $(A-\lambda^\prime I_n)u \in C_2$.
But then
\[
(A-\lambda I_n)u=(A-\lambda^\prime I_n)u+(\lambda^\prime -\lambda)u \in {\rm relint}\, C_2.
\]
Clearly,
we also have $(A-\lambda I_n)u \in C_1$.
This means that $C_1 \bigcap {\rm relint}\, C_2 \ne \emptyset$.
But $C_1$ and $C_2$ are both faces of $K$,
hence $C_2 \subseteq C_1$.
The proof is complete.\hfill$\Box$\\

By specializing Theorem 4.19 to the nonnegative matrix case,
we readily obtain the following:

\begin{coro}
Let $P$ be an $n \times n$ nonnegative matrix.
Let $I$ denote the union of all classes of $P$ that have access to some basic class.
Then for any $\lambda$,
$r < \lambda < \rho(P)$,
where $r$ denotes the largest real eigenvalue of $P$ less than $\rho(P)$
(and equals $-\infty$ if there is no such eigenvalue),
we have $\Phi((P-\lambda I_n)\IR^n_+ \bigcap \IR^n_+)=F_I$.
\end{coro}

We would like to mention that there is a direct proof of Corollary 4.20 that makes use of the trace down method and the nonnegative matrix case of Theorem 4.15.

Corollary 4.20,
in turn,
yields the following:

\begin{coro}
Let $P$ be an $n \times n$ nonnegative matrix,
and let $0 \ne b \in \IR^n_+$.
For $\lambda < \rho(P)$,
sufficiently close to $\rho(P)$,
if the equation $(P-\lambda I_n)x=b$,
$x \ge 0$ is solvable,
then for any class $\alpha$ of $P$ for which $\alpha \bigcap {\rm supp}(b) \ne \emptyset$,
we have $\alpha >= \beta$ for some basic class $\beta$ of $P$.
\end{coro}

\section{Collatz-Wielandt numbers and local spectral radii}

For any $A\in \pi (K)$ and $0 \ne x \in K$,
it is known that the local spectral radius and the lower and upper Collatz-Wielandt numbers are related by:
$$r_A(x) \leq \rho_x(A) \leq R_A(x)$$
(see,
for instance,
[T--W,
Theorem 2.4(i)]).
Clearly,
$r_A(x)=R_A(x)$ if and only if $x$ is an eigenvector of $A$.
We are going to characterize when $\rho_x (A) < R_A (x)$ with $R_A(x)$ finite,
and when $\rho_x (A)=R_A(x)$.

\begin{remark}
{\rm
Let $A \in \pi (K)$,
and let $0 \neq x\in K$.
Then the face $\Phi(x)$ is $A$-invariant if and only if $R_A (x)$ is finite.}
\end{remark}

The preceding remark is obvious,
because $\Phi (x)$ is $A$-invariant if and only if $Ax ~^K\!\!\!\leq \sigma x$ for some $\sigma>0$.

\begin{theorem}
Let $A \in \pi (K)$,
and let $0 \neq x\in K$ be such that $R_A(x)$ is finite.
Let $b$ denote the vector $(R_A(x) I_n-A)x$.
Then

{\rm (i)}
$b$ belongs to the relative boundary of $\Phi(x)$.

{\rm (ii)}
The inequality $\rho_x(A)< R_A(x)$ holds if and only if $\Phi (x)$ is the smallest A-invariant face containing b.
\end{theorem}

\noindent
{\it Proof.}
First,
note that when $R_A(x)=0$,
we have $Ax=0$,
$b=0$ and $\rho_x (A)=0$.
In this case,
(i) and (ii) clearly hold.
Hereafter,
we assume that $R_A (x)>0$.

(i)
By definition of $b$,
we have $0 ~^K\!\!\!\le b ~^K\!\!\!\leq R_A(x) x$,
so $b \in \Phi(x)$.
Indeed,
$b$ lies on the relative boundary of $\Phi (x)$,
because for any $\varepsilon >0$,
we have,
$b-\varepsilon x=((R_A(x) -\varepsilon) I_n -A)x \notin K$,
in view of the definition of $R_A (x)$.

(ii)
``Only if" part:
By the local Perron-Schaefer condition on $A$ at $x$,
there is a generalized eigenvector $y$ of $A$ corresponding to $\rho_x (A)$ that appears in the representation of $x$ as a sum of generalized eigenvectors of $A$.
Since $R_A (x)> \rho_x (A)$,
$(R_A (x) I_n -A)y$ is nonzero and so is also a generalized eigenvector of $A$ corresponding to $\rho _x (A)$ that appears in the corresponding representation for $b$.
So we have $\rho_b (A)=\rho_x (A)< R_A(x)$,
and by Remark 3.8 it follows that we have $x=(R_A (x))^{-1} \sum^{\infty}_{j=0} (R_A (x)^{-1} A)^j b$,
and hence $x \in \Phi (\hat{b})$,
where $\hat{b}=(I_n +A)^{n-1} b$.
But we also have $\hat{b}\in \Phi(x)$,
as $b\in \Phi (x)$ and $\Phi (x)$ is $A$-invariant (according to Remark 5.1).
Thus,
we have $\Phi(x)=\Phi (\hat{b})$ and by Lemma 2.1 $\Phi (x)$ is the smallest $A$-invariant face containing $b$.

``If$\,$" part:
If $\Phi(x)$ is the smallest $A$-invariant face containing $b$,
then by Lemma 2.1 $\Phi(x)=\Phi(\hat{b})$,
where $\hat{b}$ has the same meaning as above.
Since $(R_A (x) I_n -A)x =b$ and $b \in K$,
by Theorem 3.1 we have,
$R_A(x)> \rho_b (A)=\rho_{\hat{b}} (A)=\rho_x (A)$,
where the first equality holds by Lemma 2.1 and the second equality holds as
$\Phi (\hat{b})=\Phi(x)$.
This completes the proof.\hfill$\Box$\\

It is clear that by Remark 5.1 and Theorem 5.2 we have the following:

\begin{coro}
Let $A \in \pi (K)$ and let $0 \neq x\in K$.
Then $R_A (x)=\rho_x (A)$ if and only if $\Phi (x)$ is an $A$-invariant face,
but is not the smallest $A$-invariant face containing the vector $R_A(x) x -Ax$.
\end{coro}

Below we give some more explicit characterizations for when $R_A(x)=\rho_x (A)$.

\begin{theorem}
Let $A \in \pi (K)$,
and let $0 \neq x \in K$.
The following conditions are equivalent:
\vspace*{-2mm}
\leftmargini=6mm
\begin{enumerate} \itemsep=-1pt
\item[] {\rm (a)}
$R_A (x) =\rho_x (A)$.
\item[] {\rm (b)}
$(\rho_x (A) I_n -A)x \in K$.
\item[] {\rm (c)}
x can be written as $x_1+x_2$,
where $x_1, x_2 \in K$ such that $x_1$ is an eigenvector of A corresponding to $\rho_x(A)$
and $x_2$ satisfies $\rho_{x_2} (A) < \rho_x (A)$ and $R_A (x_2) \leq \rho_x (A)$.
\end{enumerate}
\end{theorem}

\noindent
{\it Proof.}
We always have the inequality $\rho_x(A) \leq R_A (x)$.
So by the definition of the upper Collatz-Wielandt number,
the equivalence of (a) and (b) follows.

(b)$\Longrightarrow$(c):
Let $b$ denote the vector $(\rho_x(A) I_n -A)x$.
If $b$ is the zero vector,
we are done.
So assume that $b \neq 0$.
By Theorem 3.1 we have $\rho_x (A)> \rho_b (A)$,
and also $(\rho_x (A) I_n -A)x^0=b$,
where $x^0=\rho_x (A)^{-1} \sum^{\infty}_{j=0} (\rho_x (A)^{-1} A)^j b \in K$.
Furthermore,
$x-x^0$ is either the zero vector or is a distinguished eigenvector of $A$ corresponding to $\rho_x(A)$.
Indeed,
the latter must happen,
as $\rho_x (A) > \rho_b (A) =\rho_{x^0} (A)$,
where the equality holds by Remark 3.8.
Since $(\rho_x (A) I_n -A)x^0 =b \geq^K 0$,
we also have $\rho_x(A) \geq R_A (x^0)$.
Set $x_1=x-x^0$ and $x_2=x^0$.
Then $x=x_1 +x_2$ is the desired decomposition for $x$.

(c)$\Longrightarrow$(b):
Straightforward.\hfill $\Box$\\

\begin{coro}
Let $A \in \pi (K)$.
If $\rho (A)\in \sum_1$,
then $\nu_{\rho(A)}(A)=1$.
\end{coro}

\noindent
{\it Proof.}
If $\rho (A) \in \sum_1$,
then there exists $x\in {\rm int}\, K$ such that $(\rho(A)I_n-A)x \in K$.
Since $x \in \,$int$\, K$,
by [T--S 2,
Lemma 4.3] we have $\rho_x (A)=\rho (A)$ and ord$_A (x) =\nu_{\rho(A)}(A)$.
By Theorem 5.4 we can write $x$ as $x_1+x_2$,
where $x_1$ is an eigenvector of $A$ corresponding to $\rho(A)$,
and $\rho_{x_2} < \rho (A)$;
hence ord$_A (x)=1$.
So we have $\nu_{\rho(A)}(A)=1$.\hfill $\Box$\\

According to [Scha,
Chapter 1,
Proposition 2.8],
if $P$ is a nonnegative matrix for which there exists a positive vector $z$ such that $P^T z \leq \rho(P)z$,
then $\rho (P)$ is a simple pole of the resolvent of $P$.
Clearly this observation also follows from our preceding corollary.

Using an argument given in the proof of [T--W,
Theorem 5.2] one readily obtains the following related result.

\begin{remark}
{\rm
Let $A\in \pi (K)$.
Let $x \in K$ with $\rho_x (A)>0$.
Then $x$ can be written as $x_1+x_2$,
where $x_1$ is an eigenvector of $A$ corresponding to $\rho_x (A)$ and $\rho_{x_2} (A) <\rho_x (A)$ if and only if $\lim_{k\to \infty} (A /\rho_x (A))^k x$ exists.}
\end{remark}

The values of the greatest lower bound or the least upper bound of the Collatz-Wielandt sets associated with $A$
($\in \pi (K)$) are known.
Specifically,
we have $\sup \Omega (A)= \inf \sum_1 (A) =\rho(A)$,
$\inf \sum(A)$ is equal to the least distinguished eigenvalue of $A$ for $K$,
and $\sup \Omega_1 (A) =\inf \sum(A^T)$ and hence is equal to the least distinguished eigenvalue of $A^T$ for $K^*$
(see [T--W,
Theorems 3.1,
3.2 and 3.3]).
(See [Fri] for an extension of these results to the settings of a Banach space or $C^*$ algebra.)
It is clear that we always have $\sup \Omega \in \Omega$ and $\inf \sum \in \sum$.
However,
in general,
$\sup \Omega_1 \notin \Omega_1$;
but when $K$ is polyhedral,
we always have $\sup \Omega_1 \in \Omega_1$ (see [Tam 1,
Example 5.5 and Corollary 5.2]).
Also,
in general,
inf$~\sum_1 \notin \sum_1$,
not even in the nonnegative matrix case (see [T--W,
Theorem 5.2] and [Tam 1,
Corollary 5.3]).

It is easy to show that $\rho(A)\in \sum_1$ if and only if there exists $x \in {\rm int}\,K$ such that $R_A(x)=\rho (A)$.
(Similarly,
we also have,
$\sup \Omega_1 \in \Omega_1$ if and only if there exists $x \in {\rm int}\,K$ such that $r_A (x)=\sup \Omega_1$.)
But such observation hardly tells anything new.
Below we give a concrete characterization for when $\inf\sum_1 \in \sum_1$.

\begin{theorem}
Let $A\in \pi (K)$ with $\rho(A)>0$.
Let $C$ denote the set $\{ x \in K: \rho_x(A) <\rho(A) \}$.
Then $\rho(A) \in \sum_1$ if and only if $\Phi (({\cal N}(\rho (A) I_n -A) \bigcap K) \bigcup C)=K$.
\end{theorem}

\noindent
{\it Proof.}
``If$\,$" part:
If $C=\{ 0 \}$,
our condition becomes $\Phi ({\cal N}(\rho (A) I_n -A) \bigcap K)=K$.
Then there exists an eigenvector $x \in {\rm int}\, K$ corresponding to $\rho(A)$,
and it is clear that we have $\rho(A) \in \sum_1$.
Hereafter,
we assume that $C \neq \{ 0\}$.

By Remark 3.5 $C$ is an $A$-invariant face of $K$.
By [T--W,
Theorem 3.1] we have $\inf \sum_1 (A|_{{\rm span}\,C})=\rho (A|_{{\rm span}\,C})$.
By definition of $C$ and an application of the Perron-Frobenius theorem to $A|_{{\rm span}\,C}$,
we also have $\rho (A) > \rho (A|_{{\rm span}\,C})$.
(In fact,
$\rho (A|_{{\rm span}C})$ is equal to the largest distinguished eigenvalue of $A$ which is less than $\rho(A)$.)
So we can choose a vector $x_2$ from relint $C$ that satisfies $\rho(A)>R_A(x_2) \geq \rho(A|_{{\rm span}C})$.
Now choose any vector $x_1$ from relint$({\cal N} (\rho (A) I_n -A)\bigcap K)$,
and let $x=x_1+x_2$.
In view of the condition $\Phi(({\cal N}(\rho (A) I_n -A) \bigcap K) \bigcup C)=K$,
we have $x \in {\rm int} \, K$.
By our choices of $x_1$ and $x_2$,
we also have,
$(\rho (A) I_n -A)x=(\rho (A) I_n -A)x_2 \in K$.
So we have $\rho(A) \in \sum_1$.

``Only if$\,$" part:
Let $x \in \,$int$\,K$ be such that $(\rho (A) I_n -A)x \in K$.
Since $x \in \,$int$\,K$,
we have $\rho_x (A)=\rho (A)$.
By Theorem 5.4 $x$ can be written in the form $x_1+x_2$,
where $x_1, x_2 \in K$,
$x_1$ is an eigenvector of $A$ corresponding to $\rho (A)$,
and $x_2$ satisfies $\rho_{x_2} (A)< \rho(A)$ (and $R_A (x_2) \leq \rho (A)$).
But $x_1+x_2 \in \Phi(({\cal N}(\rho (A)I_n-A)\bigcap K) \bigcup C)$,
so we have $\Phi(({\cal N}(\rho (A)I_n-A)\bigcap K) \bigcup C)=K$.\\
\hspace*{13.4cm}$\Box$

\begin{remark}
{\rm
In case $\rho(A)$ is the only distinguished eigenvalue of $A$ for $K$,
the equivalent condition given in Theorem 5.7 for $\rho (A) \in \sum_1$ reduces to
``$A$ has an eigenvector in int $K$ (corresponding to $\rho (A)$)".}
\end{remark}

For the nilpotent case we have the following obvious result.

\begin{remark}
{\rm
If $A \in \pi (K)$ is nilpotent,
then $0 \in \sum_1$ if and only if $A=0$.}
\end{remark}

Now we rederive the corresponding known result for the nonnegative matrix case
(see [Sch 3,
Theorem 5.1] or [T--W,
Theorem 5.2]):

\begin{theorem}
Let P be an $n \! \times \! n$ nonnegative matrix.
A necessary and sufficient condition for $\rho(P) \in \sum_1$ is that every basic class of P is final.
\end{theorem}

\noindent
{\it Proof.}
If $P$ is nilpotent,
then by Remark 5.9 we have $\rho (P) \in \sum_1$ if and only if $P$ is the zero matrix.
In this case,
each class of $P$ is a singleton and is also basic.
So it is clear that the condition
``every basic class is final" is equivalent to $P$ being the zero matrix.
This proves our assertion for the nilpotent case.

Suppose that $P$ is non-nilpotent.
By Theorem 5.7 and in its notation
(but with $A$ and $K$ replaced by $P$ and $\IR^n_+$ respectively)
we have
\[
\rho(P) \in {\textstyle\sum_1} \mbox{~~if and only if~~} \Phi(({\cal N}(\rho (P)I_n-P){\textstyle\bigcap} \IR^n_+) {\textstyle\bigcup} C)=\IR^n_+.
\]
Note that $\Phi({\cal N}(\rho (P)I_n-P)\bigcap \IR^n_+)$ and $C$ are both $P$-invariant faces of $\IR^n_+$;
so by Theorem 2.2 they can be written as $F_{I_1}$ and $F_{I_2}$ respectively,
where $I_1$ and $I_2$ are initial subsets for $P$.

>From the proof of Corollary 4.2,
$I_1$ is equal to the union of all classes that have access to a distinguished basic class.

By the beginning part of the proof of Corollary 3.3,
for any vector $x \in \IR_+^n$,
we have $\rho_x (P)=\max \{\rho(P_{\alpha\alpha}):~\alpha$ has access to supp$(x)\}$.
In view of the definition of $C$,
for any $x \in K$,
we have $x \in C$ if and only if supp$(x)$ has no access from a basic class if and only if supp$(x)$ is included in the union of all classes that have no access from a basic class.
So it is clear that $I_2$ is equal to the union of all classes that have no access from a basic class.

Now we have
\begin{eqnarray*}
\Phi(({\cal N}(\rho (P)I_n-P){\textstyle\bigcap} \IR^n_+) {\textstyle\bigcup} C)
&=& \Phi({\cal N}(\rho (P)I_n-P){\textstyle\bigcap} \IR^n_+) {\textstyle\bigvee} C \\
&=& F_{I_1} {\textstyle\bigvee} F_{I_2} \\
&=& F_{I_1 \cup I_2}~.
\end{eqnarray*}
So $\rho(P) \in \sum_1$ if and only if $F_{I_1 \cup I_2}=\IR^n_+$ if and only if $I_1 \bigcup I_2 =\langle n \rangle$.

Suppose that each basic class of $P$ is final.
Consider any final class $\alpha$ of $P$.
If $\alpha$ is basic,
then $\alpha$ must be distinguished;
otherwise,
we would have a basic class which is not a final class.
Hence,
$\alpha$ is included in $I_1$.
If $\alpha$ is non-basic,
then $\alpha$ cannot have access from a basic class and so $\alpha$ is included in $I_2$.
This shows that each final class of $P$ is included in $I_1 \bigcup I_2$.
But $I_1 \bigcup I_2$ is an initial subset for $P$,
it follows that $I_1 \bigcup I_2 =\langle n \rangle$ and hence $\rho(P) \in \sum_1$.

Conversely,
suppose $\rho(P) \in \sum_1$,
or equivalently,
$I_1 \bigcup I_2 =\langle n \rangle$.
Consider any basic class $\alpha$ of $P$.
If $\alpha \bigcap I_2 \neq \emptyset$,
then since $I_2$ is an initial subset for $P$,
we would obtain $\alpha \subseteq I_2$,
a contradiction.
So $\alpha \subseteq I_1$;
hence $\alpha$ is a distinguished basic class,
and moreover $\alpha$ is final in the initial collection of classes corresponding to $I_1$.
Suppose that $\alpha$ is not a final class of $P$.
Then $\alpha >\!\!-\, \beta$ for some class $\beta$.
By definition of $I_2$,
$\beta$ must be disjoint from $I_2$.
So $\beta$ is included in $I_1$.
But then $\alpha$ is not final in the initial collection of classes corresponding to $I_1$,
which is a contradiction.
Therefore,
$\alpha$ must be a final class of $P$.
The proof is complete.\hfill $\Box$\\

We take a digression and return to the question of when $K \bigcap (A-\rho(A) I_n)K$ $= \{0\}$.
We have the following result,
which contains Corollary 4.10 as well as [Tam 1,
Corollary 4.3].

\begin{theorem}
Let $A\in \pi (K)$.
Consider the following conditions:

{\rm (a)}
$\rho(A) \in \sum_1 (A^T)$.

{\rm (b)}
${\cal N}((\rho(A) I_n -A)^n) \bigcap K ={\cal N}(\rho(A) I_n-A)\bigcap K$,
and $A$ has no eigenvectors in $\Phi({\cal N}(\rho(A) I_n -A)\bigcap K)$ corresponding to an eigenvalue other than $\rho(A)$.

{\rm (c)}
$K \bigcap (A-\rho(A)I_n)K=\{0\}$.\\
We always have {\rm (a)}$\Longrightarrow${\rm (b)}$\Longrightarrow${\rm (c)}.
When $K$ is polyhedral,
conditions {\rm (a)},
{\rm (b)} and {\rm (c)} are equivalent.
\end{theorem}

\noindent
{\it Proof}.
(a)$\Longrightarrow$(b):
Since $\rho(A) \in \sum_1 (A^T)$,
by Theorem 5.7 we have $K^*=\Phi(({\cal N}(\rho(A)I_n-A^T)\bigcap K^*)\bigcup F)$,
where we use $F$ to denote the face $\{z \in K^*: \rho_z (A) < \rho(A) \}$ of $K^*$.
In this case,
by Corollary 5.5 we also have $\nu_{\rho(A)}(A)=1$;
hence,
the first part of condition (b) holds.
Assume to the contrary that $A$ has an eigenvector in $\Phi({\cal N}(\rho(A)I_n -A)\bigcap K)$ corresponding to a (distinguished) eigenvalue,
say $\lambda$,
other than $\rho(A)$.
For simplicity,
we denote by $G$ the face $\Phi({\cal N}(\lambda I_n -A)\bigcap K)\bigcap \Phi({\cal N}(\rho(A)I_n -A)\bigcap K)$.
[Note that we need not have $G=\Phi({\cal N}(\lambda I_n-A)\bigcap K)$.]
Then we have $d_K(\Phi({\cal N}(\rho(A)I_n-A)\bigcap K))\subseteq d_K(G)$.
Clealy,
every vector in $F$ is orthogonal to ${\cal N}(\rho(A)I_n -A)\bigcap K$.
So we must have $F \subseteq d_K(\Phi({\cal N}(\rho(A)I_n-A)\bigcap K))$.
Similarly,
we also have $\Phi({\cal N}(\rho(A)I_n -A^T)\bigcap K^*)\subseteq d_K(G)$.
Hence,
we have
\[
d_K(G) \supseteq \Phi([{\cal N}(\rho(A)I_n -A^T){\textstyle \bigcap} K^*]{\textstyle \bigcup} F)=K^*,
\]
which is a contradiction,
as $G$ is a nonzero face of $K$.

(b)$\Longrightarrow$(c):
Follows form Corollary 4.10.

If $K$ is polyhedral,
the implication (c)$\Longrightarrow$(a) follows from [Tam 1,
Corollary 4.3].
Then conditions (a)--(c) are equivalent.\hfill$\Box$\\

When $K$ is non-polyhedral,
the missing implications in theorem 5.11 all do not hold.
For instance,
the matrix $A$ considered in Example 4.9 satisfies condition (b),
but it does not satisfies condition (a),
in view of Corollary 5.5,
as $\nu_{\rho(A)}(A) \ne 1$.
So we have (b)$\narrow$(a).
Example 4.11 can also be used to illustrate (c)$\narrow$(b).

For the question of when $r_A(x)=\rho_x(A)$,
we have two partial results.
First,
it is straightforward to show the following:

\begin{remark}
{\rm
Let $A\in \pi (K)$,
and let $0 \ne x \in K$.
Then $r_A (x)=\rho_x (A)$ if and only if $(A-\rho_x (A) I_n)x \in K$.}
\end{remark}

\begin{theorem}
Let $A\in \pi (K)$,
and let $0 \ne x \in K$.
Then we have ord$_A(x)=1$ and $r_A(x)=\rho_x (A)$ if and only if $x$ can be written as $x_1 -x_2$,
where $x_1$,
$x_2 \in K$ such that $x_1$ is an eigenvector of $A$ corresponding to $\rho_x (A)$ and $x_2$ satisfies $\rho_{x_2} (A) < \rho_x (A)$ and $R_A (x_2) \leq \rho_x (A)$.
\end{theorem}

\noindent
{\it Proof.}
It is straightforward to verify the ``if" part,
using Remark 5.12.

To prove the ``only if" part,
suppose that we have ord$_A (x)=1$ and $r_A (x)=\rho_x (A)$.
By Remark 5.12,
the vector $(A-\rho_x (A) I_n)x$,
which we denote by $b$,
belongs to $K$.
Since ord$_A (x)=1$,
by using the local Perron-Schaefer conditions on $A$ at $x$ and $b$ respectively,
we readily obtain $\rho_b (A) < \rho_x (A)$.
Let $x_2$ denote the vector $\sum^{\infty}_{k=0} (\rho_x (A))^{-k-1} A^k b$.
Also,
let $x_1=x+x_2$.
By the proof of the ``only if" part of Theorem 4.1
(with $\rho_x (A),~x_1,~x_2$ in place of $\lambda, ~w, ~x_0$ respectively),
we find that $x_1, ~x_2$ are both vectors of $K$ such that $x_1$ is an eigenvector of $A$ corresponding to $\rho_x (A)$ and $\rho_{x_2} (A)=\rho_b (A)< \rho_x (A)$.
Since $\rho_x (A) x_2 -Ax_2 =(A-\rho_x (A) I_n)x \in K$,
we also have $\rho_x (A) \geq R_A(x_2)$.
Thus,
$x=x_1 -x_2$ is the desired decomposition for $x$.\hfill$\Box$\\

In case $x \in K$ satisfies ord$_A(x) \geq 2$,
we do not know when $r_A(x)=\rho_x (A)$ holds.

The question of when $\sup \Omega_1 \in \Omega_1$ seems to be more subtle than that of when $\inf \sum_1 \in \sum_1$.
It is easy to see that when $\rho (A) (= \inf\sum_1) \in \sum_1$ and $x \in {\rm int}\,K$ satisfies $\rho (A) x \geq^K Ax$,
necessarily we have $\rho_x (A) =R_A (x)$.
This explains why Theorem 5.4 has been useful in proving Theorem 5.7.
In contrast,
when $\sup\Omega_1 \in \Omega_1$ and $x \in {\rm int}\,K$ satisfies $Ax \geq^K (\sup \Omega_1)x$,
we only have $\sup \Omega_1 =r_A (x) \leq \rho_x (A) =\rho (A)$.
In particular,
if $\sup\Omega_1$,
which is the least distinguished eigenvalue of $A^T$ for $K^*$,
is less than $\rho (A)$,
then we cannot expect that a solution for the question of when $r_A(x)=\rho_x (A)$ is of help in answering when $\sup \Omega_1 \in \Omega_1$.
Indeed,
in this case we are faced with equation (1.2)
(with $\lambda=\sup \Omega_1$ and $b=Ax -(\sup\Omega_1)x$)
for the case when $\lambda< \rho_b (A)$,
which is the case we know not so well.

\section{Alternating sequences}

Let $A \in {\cal M}_n (\IC)$ and let $x \in \IC^n$.
Following [H--R--S],
we call the sequence $x, Ax, \ldots, A^k x$ an {\it alternating sequence for A of length k} if $(-1)^k A^k x \geq 0$ and $0 \neq (-1)^r A^r x \geq 0$ for $r=0,\ldots, k-1$.
Here $\geq$ denotes the usual componentwise partial ordering of $\IR^n$.
The infinite sequence $x, Ax, A^2x, \ldots$ is said to be an {\it infinite alternating sequence for A} if $0 \neq (-1)^r A^r x \geq 0$ for $r=0,1,2,\ldots$.

Recall that an $n \! \times \! n$ real matrix is called a $Z$-matrix if it is of the form $\lambda I_n -P$,
where $P$ is nonnegative.
In [H--R--S,
Corollary 3.5]
the following characterization of $M$-matrices among $Z$-matrices is given:

{\it Let $A$ be a $Z$-matrix.
Then $A$ is an $M$-matrix if and only if every alternating sequence for $A$ is of finite length.}

Using the local Perron-Schaefer conditions on a nonnegative matrix,
we can readily explain why the above result is true.
Indeed,
we can extend the result to the setting of a cone-preserving map.

\begin{theorem}
Let $A \in \pi (K)$,
let $0 \neq x\in K$,
and let $x=x_1 +\cdots +x_k$ be the representation of $x$ as a sum of generalized eigenvectors of $A$,
where $\lambda_1, \ldots, \lambda_k$ are the corresponding distinct eigenvalues.
Let $\Gamma$ denote the set $\{ j\in \langle k\rangle :| \lambda_j |=\rho_x (A)$ and $\lambda_j \neq \rho_x(A) \}$.
Let $m$ be a positive integer and suppose that $(A-\rho_x (A) I_n)^m x \in K$ and $0\neq (A-\rho_x (A) I_n)^j x \in K$ for $j=0,\ldots, m-1$.
If $\Gamma = \emptyset$,
then $m \leq {\rm ord}_A (x)$.
If $\Gamma \neq \emptyset$,
then $m\leq {\rm ord}_A (x) -\max_{j \in \Gamma} {\rm ord}_A (x_j)$.
\end{theorem}

\noindent
{\it Proof.}
By the local Perron-Schaefer condition on $A$ at $x$ there is an index $j$ such that $\lambda_j=\rho_x (A)$ and ord$_A (x_j) ={\rm ord}_A (x)$.
Let $y$ denote the vector $(A-\rho_x (A) I_n)^{{\rm ord}_A (x)} x$.
If $m > {\rm ord}_A (x)$,
then,
by our hypothesis,
$y$ is a nonzero vector of $K$.
Clearly,
we do not have a generalized eigenvector corresponding to $\rho_x(A)$ that appears in the representation of $y$ as a sum of generalized eigenvectors of $A$.
So by the local Perron-Schaefer condition on $A$ at $y$,
we have $\rho_y (A)< \rho_x (A)$.
But we also have $\rho_y (A) \geq r_A (y)$;
hence $(A-\rho_x (A) I_n)^{{\rm ord}_A (x)+1} x=(A-\rho_x (A) I_n)y \notin K$,
in contradiction to our hypothesis.
This proves that we always have $m \leq {\rm ord}_A (x)$.

It remains to consider the case when $\Gamma \neq \emptyset$.
Denote by $t$ the value of $\max_{j\in \Gamma}$ ${\rm ord}_A (x_j)$.
According to the local Perron-Schaefer condition on $A$ at $x$,
we have $t \leq {\rm ord}_A (x)$.
Let $w$ denote the vector $(A-\rho_x (A) I_n)^{{\rm ord}_A (x)-t+1} (x)$,
and consider its representation as a sum of generalized eigenvectors of $A$.
Note that in the representation there is at least a generalized eigenvector of order $t$ corresponding to an eigenvalue different from $\rho_x (A)$ but with modulus $\rho_x (A)$,
and also that if $t>1$ then the order of the generalized eigenvector corresponding to $\rho_x (A)$ that appears in the representation is $t-1$,
and if $t=1$ then in the representation there does not exist a generalized eigenvector coresponding to $\rho_x(A)$.
Hence,
we have,
$\rho_w (A) =\rho_x (A)$ and the local Perron-Schaefer condition on $A$ at $w$ is not satisfied.
It follows that $w \notin K$.
So by our hypothesis,
we have,
$m<{\rm ord}_A (x) -t+1$,
i.e. $m \leq {\rm ord}_A (x) -t$,
which is the desired inequality.\hfill $\Box$\\

\begin{coro}
Let $A\in \pi (K)$,
let $x \in K$,
and let $\lambda$ be a real number.
In order that we have $0 \neq (A-\lambda I_n)^j x \in K$ for all positive integers $j$,
it is necessary that $\lambda < \rho_x (A)$.
\end{coro}

\noindent
{\it Proof.}
Assume that $0 \neq (A-\lambda I_n)^j x \in K$ for $j=0, 1, 2, \ldots$.
By the conditions $(A-\lambda I_n)x \in K$ and $0 \neq x \in K$,
clearly we have $\lambda \leq r_A (x) \leq \rho_x (A)$.
If $\lambda=\rho_x (A)$,
then by Theorem 6.1 we would have ord$_A (x) \geq m$ for each positive integer $m$,
which is impossible.
So we must have $\lambda < \rho_x (A)$.\hfill $\Box$\\

Our next result is an extension of [H--R--S,
Corollary 3.5]
(which was mentioned at the beginning of this section):

\begin{coro}
Let $A\in \pi (K)$,
and let $\lambda$ be a real number.
Then $\lambda<\rho(A)$ if and only if there exists a vector $x \in K$ such that $0 \neq (A-\lambda I_n)^j x \in K$ for all positive integers $j$.
\end{coro}

\noindent
{\it Proof.}
``If$\,$" part:
By Corollary 6.2 we have $\lambda < \rho_x (A) \leq \rho(A)$.

``Only if$\,$" part:
Suppose $\lambda <\rho(A)$.
Let $x \in K$ be an eigenvector of $A$ corresponding to $\rho(A)$.
Then $(A-\lambda I_n)^j x=(\rho (A) -\lambda)^j x$ is a nonzero vector of $K$ for all positive integers $j$.\hfill $\Box$\\

In [H--R--S,
Theorem 3.4(ii)] it is proved that if $A$ is an $M$-matrix,
then the index of $A$ is equal to the maximal length of an alternating sequence for $A$.
Making use of Theorem 6.1,
we readily obtain the following partial extension:

\begin{coro}
Let $A\in \pi (K)$,
let $x \in K$,
and let $m$ be a positive integer.
If $(A-\rho (A) I_n)^m x \in K$ and $0\neq (A-\rho (A) I_n)^j x \in K$ for $j=0, \ldots, m-1$,
then $\rho_x (A)=\rho (A)$ and $m \leq {\rm ord}_A (x)\leq \nu_{\rho(A)} (A)$.
\end{coro}

According to [Tam 1,
Theorem 7.5],
if $K$ is a polyhedral cone,
then for any $A \in \pi (K)$,
there exists a vector $x \in K$ such that $(A-\rho(A) I_n)^{\nu} x=0$ and $0 \neq (A-\rho(A) I_n)^j x \in K$ for $j=1,\ldots, \nu -1$,
where $\nu=\nu_{\rho(A)} (A)$.
So in the polyhedral case,
we have a full extension of [H--R--S,
Theorem 3.4(ii)].

In the nonpolyhedral case,
the other extreme can happen.
Example 4.11 can be used to show that for any odd integer $n \ge 3$,
there exists a proper cone $K$ and a matrix $A \in \pi(K)$ such that $\nu_{\rho(A)}(A)=n$,
and for any $0\ne x \in K$,
we have $(A-\rho(A)I_n)x \notin K$,
except when $x$ is an eigenvector of $A$ corresponding to $\rho(A)$.\\

{\it
Thanks are due to Ludwig Elsner for observing,
after the first author's talk at Oberwolfach,
that $(tI_n -A)^{-1}$ is negative if $A$ is irreducible nonnegative and $t$ is slightly less than $\rho(A)$,
which has aroused further work in the later parts of Section 4.}

\section*{\bf Appendix A.  A proof for {\boldmath$\rho_x(A)=\lim_{m\to \infty} \| A^m x \|^{1/m}$} }

{\bf Proof.}
To show this,
we may assume that $W_x =\IC^n$.
Let $x=x_1 + \cdots + x_k$ be the representation of $x$ as a sum of generalized eigenvectors of $A$ corresponding to the distinct eigenvalues $\lambda_1, \ldots , \lambda_k$ respectively.
For each $i=1, \ldots, k$,
let $n_i$ denote the order of the generalized eigenvector $x_i$.
Let $J_q (\lambda)$ denote the $q \! \times \! q$ upper triangular elementary Jordan matrix corresponding to $\lambda$.
Then the vectors $x_j, Ax_j,\ldots, A^{n_j-1}x_j$,
$j=1,\ldots,k$,
constitute a basis for $\IC^n$,
and we can find a nonsingular matrix $P$ such that $P^{-1}AP=J_{n_1} (\lambda_1) \bigoplus \cdots \bigoplus J_{n_k} (\lambda_k)$,
and $P^{-1} x$ is the vector of $\IC^n$ with 1's at its $n_1$th,
$(n_1+n_2)$th,
$\ldots$,
and $(n_1+\cdots+n_k)$th components and with 0's elsewhere.
Since any two norms on $\IC^n$ are equivalent,
the existence of the limit $\lim _{m\to \infty} \|A^m x \|^{1/m}$ and also its value are independent of the choice of the norm of $\IC^n$.
Here we choose the norm on $\IC^n$ given by $\|y\| =\| P^{-1} y \|_1$,
where $\| \cdot \|_1$ denotes the $l_1$-norm.
For all positive integers $m$,
we have $\|A^m x \|=\|J^m P^{-1} x\|_1$.
Using the fact that $\rho_x(A)=\max _{1\leq j\leq k} | \lambda_j |$,
and after a little calculation,
we readily see that for all positive integers $m>n$,
we have
\[
\rho_x (A)^m \leq \| J^m P^{-1} x \|_1 \leq n\rho_x (A)^m \frac{m!}{n!(m-n)!}.
\]
It follows that we have $\lim _{m\to \infty} \| A^m x \|^{1/m}=\rho_x (A)$.\hfill $\Box$

\section*{\bf Appendix B.  Two proofs for Theorem 3.1,
{\boldmath $(a) \Rightarrow (b)$}}

{\bf First Proof.}
When condition $(a)$ is fulfilled,
it is clear that $b\in (\lambda I_n -A)K \bigcap K$.
It is also clear that $A$ leaves invariant the nonzero,
closed pointed cone $(\lambda I_n -A)K \bigcap K$.
Let $\rho$ denote the spectral radius of the restriction of $A$ to span$((\lambda I_n -A)K \bigcap K)$.

We contend that $\rho < \lambda$.
Assume to the contrary that the reverse inequality holds.
By the Perron-Frobenius theorem,
there exists a nonzero vector $y \in (\lambda I_n -A)K \bigcap K$ such that $(\rho I_n -A)y=0$.
Then $(\lambda I_n -A)x=y$ for some vector $x \in K$.
Note that the vectors $x, \, y$ are linearly independent;
if not,
we have $y=(\lambda I_n -A)x=(\lambda -\rho )x ~^K\!\!\!\leq 0$,
which is a contradiction.
Let $C$ denote the 2-dimensional cone $K \bigcap {\rm span} \{x, y\}$.
Clearly $A$ leaves $C$ invariant.
If $\rho =\lambda$,
then the spectrum of $A|_{{\rm span}C}$ is $\{ \lambda \}$
($x$ being a generalized eigenvector of order two)
and so $\lambda$ is the only distinguished eigenvalue of $A|_{{\rm span}C}$ for $C$.
Hence,
by [Tam 1,
Theorem 5.11],
the inequality $\lambda x \geq^C Ax$ implies that $0=\lambda x-Ax=y$,
which is a contradiction.
So we have $\lambda <\rho$.
But then,
as can be readily checked by direct calculation,
$(\rho -\lambda)x+y$ is an eigenvector of $A$ corresponding to $\lambda$,
that lies in relint $C$.
Hence $\lambda =\rho (A|_{{\rm span}C})=\rho$,
which is again a contradiction.
This proves our contention.

Since $b \in (\lambda I_n -A)K \bigcap K , ~\rho_b (A) \leq \rho$.
So we have $\rho_b (A) < \lambda$.\hfill $\Box$\\

{\bf Second Proof.}
Suppose that there exists a vector $x\in K$ such that $(\lambda I_n-A)x=b$.
Multiplying both sides of the equation by $(I_n+A)^{n-1}$,
we obtain $(\lambda I_n -A)\hat{x}=\hat{b}$,
where $\hat{x}=(I_n+A)^{n-1} x,~\hat{b}=(I_n+A)^{n-1} b$,
and $\hat{x},\hat{b} \in K$.
The latter equation implies that $\lambda \in \sum_1 (A|_{W_x})$,
noting that $A|_{W_x} \in \pi (\Phi (\hat{x}))$ in view of Lemma 2.1.
By [T--W,
Theorem 3.1] and Lemma 2.1 again,
we have $\lambda \geq \inf \sum_1 (A|_{W_x})=\rho(A|_{W_x})=\rho_x(A)$.
Since $b \in W_x$,
we have $W_b \subseteq W_x$ and hence $\rho_b (A) \leq \rho_x (A)$.
So we have $\lambda \geq \rho_b(A)$.

If $\lambda > \rho_x (A)$,
we already have $\lambda > \rho_b(A)$.
Hereafter,
we assume that $\lambda=\rho_x(A)$.
Let $x=x_1+\cdots+x_k$ be the representation of $x$ as a sum of generalized eigenvectors of $A$ corresponding to distinct eigenvalues $\lambda_1,\ldots,\lambda_k$ respectively.
By the local Perron-Schaefer condition on $A$ at $x$,
we may assume that $\lambda_1=\rho_x(A)$;
then we have ord$_A(x) ={\rm ord}_A(x_1)$.
Denote this common value by $m$.

We contend that $m=1$.
Suppose that $m \geq 2$.
By [T--S 2,
Corollary 4.8]
or [Sch 2,
Theorem 5.2],
$(A-\rho_x (A) I_n)^{m-1} x_1$ is a distinguished eigenvector of $A$ corresponding to $\rho_x (A)$.
Note that the representation of $b$ as a sum of generalized eigenvectors of $A$ is
\[
(\rho_x (A) I_n -A)x_1 +(\rho_x (A) I_n -A)x_2 +\cdots+(\rho_x (A) I_n -A)x_k ,
\]
and the corresponding eigenvalues are still $\lambda_1,\ldots,\lambda_k$.
Clearly,
we have
\[
{\rm ord}_A ((\rho_x (A) I_n -A)x_1) ={\rm ord}_A (x_1) -1=m-1 \ge 1.
\]
Since $0\neq b \in K$,
by the local Perron-Schaefer condition on $A$ at $b$,
it follows that we have $\rho_b(A)=\rho_x (A)$ and ord$_A(b)=m-1$.
Applying [T--S 2,
Corollary 4.8]
to the vector $b$,
we find that $(A-\rho_x (A) I_n)^{m-2} (\rho_x (A) I_n -A)x_1=-(A-\rho_x (A) I_n)^{m-1} x_1$ is also a distinguished eigenvector of $A$ corresponding to $\rho_x(A)$.
Hence the nonzero vector $(A-\rho_x (A) I_n)^{m-1}x$ and its negative both belong to $K$,
which is a contradiction.
This proves our contention that $m=1$.

Note that now in the representation of $b$ as a sum of generalized eigenvectors of $A$,
we do not have a term which is a generalized eigenvector of $A$ corresponding to $\rho_x(A)$.
Hence,
by the local Perron-Schaefer condition on $A$ at $b$,
$\rho_b(A) \neq \rho_x(A)$.
But we always have $\rho_b(A) \leq \rho_x(A)$,
so we have $\rho_b(A) <\rho_x(A) =\lambda$.\hfill $\Box$\\

\newpage
\begin{itemize}
\item[]
\hspace{-1cm}
{\bf References}
\begin{enumerate}
\item[{[Bar]}]
G.P. Barker,
Theory of cones,
{\it Linear Algebra Appl.}
{\bf 39} (1978), 263-291.

\item[{[B--P]}]
A. Berman and R.J. Plemmons,
{\it Nonnegative Matrices in the Mathematical Sciences},
Revised reprint of the 1979 original,
Classics in Applied Mathematics,
9,
SIAM, Philadelphia, 1994.

\item[{[Car]}]
D.H. Carlson,
A note on $M$-matrix equations,
{\it SIAM J. Appl. Math.}
{\bf 11} (1963),
1027--1033.

\item[{[Fan]}]
Ky Fan,
Some matrix inequalities,
{\it Abh. Math. Sem. Univ. Hamburg} {\bf 29} (1966),
185--196.

\item[{[F--N 1]}]
K.-H. F$\ddot{\rm o}$rster and B. Nagy,
On the local spectral theory of positive operators,
{\it Oper. Theory,
Adv. Appl.}
{\bf 28} (1988),
71--81.

\item[{[F--N 2]}]
K.-H. F$\ddot{\rm o}$rster and B. Nagy,
On the Collatz-Wielandt numbers and the local spectral radius of a nonnegative operator,
{\it Linear Algebra Appl}. {\bf 120} (1989),
193--205.

\item[{[Fri]}]
S. Friedland,
Characterizations of spectral radius of positive elements on $C^*$ algebras,
{\it J. Funct. Anal.} {\bf 97} (1991),
64--70.

\item[{[Fro 1]}]
G.F. Frobenius,
$\ddot{\rm U}$ber Matrizen aus positiven Elementen,
{\it S.-B. Preuss. Akad. Wiss $($Berlin$)$} (1908),
471--476.

\item[{[Fro 2]}]
G.F. Frobenius,
$\ddot{\rm U}$ber Matrizen aus positiven Elementen,
II,
{\it S.-B. Preuss. Akad. Wiss $($Berlin$)$} (1909),
514--518.

\item[{[Fro 3]}]
G.F. Frobenius,
$\ddot{\rm U}$ber Matrizen aus nicht negativen Elementen,
{\it Sitzungsber. K$\ddot{\rm o}$n. Preuss. Akad. Wiss. Berlin},
1912,
456--477;
{\it Ges. Abh} {\bf 3},
Springer-Verlag,
1968,
546--567.

\item[{[F--S]}]
S. Friedland and H. Schneider,
The growth of powers of a nonnegative matrix,
{\it SIAM J. Algebraic Discrete Methods} {\bf 1} (1980),
185--200.

\item[{[G--L]}]
I.M. Glazman and Ju. I. Ljubi$\breve{c}$,
{\it Finite-Dimensional Linear Analysis:
A Systematic Presentation in Problem Form},
(G.P. Barker and G. Kuerti,
Transl. and Ed.),
MIT Press,
Cambridge,
Mass., 1974.

\item[{[H--R--S]}]
D. Hershkowitz,
U.G. Rothblum and H. Schneider,
Characterizations and classifications of $M$-matrices using generalized nullspaces,
{\it Linear Algebra Appl.}
{\bf 109} (1988),
59--69.

\item[{[H--S 1]}]
D. Hershkowitz and H. Schneider,
On the generalized nullspace of $M$-matrices and $Z$-matrices,
{\it Linear Algebra Appl.} {\bf 106} (1988),
5--23.

\item[{[H--S 2]}]
D. Hershkowitz and H. Schneider,
Solutions of $Z$-matrix equations,
{\it Linear Algebra Appl.}
{\bf 106} (1988), 25-38.

\item[{[J--V 1]}]
R.J. Jang and H.D. Victory,
Jr.,
Frobenius decomposition of positive compact operators,
{\it Positive Operators,
Riesz Spaces,
and Economics},
Springer Studies in Economic Theory,
Vol. 2, 195--224,
Springer Verlag,
New York,
1991.

\item[{[J--V 2]}]
R.J. Jang and H.D. Victory,
Jr.,
On nonnegative solvability of linear integral equations,
{\it Linear Algebra Appl}.
{\bf 165} (1992),
197--228.

\item[{[J--V 3]}]
R.J. Jang and H.D. Victory,
Jr.,
On the ideal structure of positive,
eventually compact linear operators on Banach lattices,
{\it Pacific J. Math.}
{\bf 157} (1993),
57--85.

\item[{[J--V 4]}]
R.J. Jang-Lewis and H.D. Victory,
Jr.,
On nonnegative solvability of linear operator equations,
{\it Integral Equations Operator Theory}.
{\bf 18} (1994),
88--108.

\item[{[Mar]}]
I. Marek,
Collatz-Wielandt numbers in general partially ordered spaces,
{\it Linear Algebra Appl.} {\bf 173} (1992),
165--180.

\item[{[Nel 1]}]
P. Nelson,
Jr.,
Subcritically for transport of multiplying particles in a slab,
{\it J. Math. Anal. Appl.}
{\bf 35} (1971),
90-104.

\item[{[Nel 2]}]
P. Nelson,
Jr.,
Positive solutions of positive linear equations,
{\it Proc. Amer. Math. Soc.}
{\bf 31} (1972),
453--457.

\item[{[Nel 3]}]
P. Nelson,
Jr.,
The structure of positive linear integral operator,
{\it J. London Math. Soc.}
{\bf 8 (2)} (1974),
711--718.

\item[{[Per 1]}]
O. Perron,
Grundlagen f$\ddot{\rm u}$r eine Theorie des Jacobischen Kettenbruchalogithmus,
{\it Math. Ann.} {\bf 63} (1907),
1--76.

\item[{[Per 2]}]
O. Perron,
Zur Theorie der $\ddot{\rm U}$ber Matrizen,
{\it Math. Ann.} {\bf 64} (1907),
248--263.

\item[{[Rot]}]
U.G. Rothblum,
Algebraic eigenspaces of non-negative matrices,
{\it Linear Algebra Appl.}
{\bf 12} (1975),
281--292.

\item[{[Scha]}]
H.H. Schaefer,
{\it Banach Lattices and Positive Operators},
Springer,
New York,
1974.

\item[{[Sch 1]}]
H. Schneider,
The elementary divisors,
associated with 0,
of a singular $M$-matrix,
{\it Proc. Edinburgh Math. Soc.},
{\bf 10} (1956),
108--122.

\item[{[Sch 2]}]
H. Schneider,
Geometric conditions for the existence of positive eigenvalues of matrices,
{\it Linear Algebra Appl.}
{\bf 38} (1981),
253--271.

\item[{[Sch 3]}]
H. Schneider,
The influence of the marked reduced graph of a nonnegative matrix on the Jordan form and on related properties:
a survey,
{\it Linear Algebra Appl.}
{\bf 84} (1986),
161--189.

\item[{[Tam 1]}]
B.S. Tam,
On the distinguished eigenvalues of a cone-preserving map,
{\it Linear Algebra Appl.}
{\bf 131} (1990),
17--37.

\item[{[Tam 2]}]
B.S. Tam,
A cone-theoretic approach to the spectral theory of positive linear operators:
the finite-dimensional case,
{\it Taiwanese J. Math.} {\bf 5} (2001),
207--277.

\item[{[Tam 3]}]
B.S. Tam,
On matrices with invariant closed,
pointed cones,
in preparation.

\item[{[T--S 1]}]
B.S. Tam and H. Schneider,
On the core of a cone-preserving map,
{\it Trans. Amer. Math. Soc.}
{\bf 343} (1994),
479--524.

\item[{[T--S 2]}]
B.S. Tam and H. Schneider,
On the invariant faces associated with a cone-preserving map,
{\it Trans. Amer. Math. Soc.} {\bf 353} (2001),
209--245.

\item[{[T--W]}]
B.S. Tam and S.F. Wu,
On the Collatz-Wielandt sets associated with a cone-preserving map,
 {\it Linear Algebra Appl.}
{\bf 125} (1989),
77--95.

\item[{[Vic 1]}]
H.D. Victory,
Jr.,
On linear integral operators with nonnegative kernels,
{\it J. Math. Anal. Appl.}
{\bf 89} (1982),
420--441.

\item[{[Vic 2]}]
H.D. Victory,
Jr.,
The structure of the algebraic eigenspace to the spectral radius of eventually compact,
nonnegative integral operators,
{\it J. Math. Anal. Appl.}
{\bf 90} (1982),
484--516.

\item[{[Vic 3]}]
H.D. Victory,
Jr.,
On nonnegative solutions to matrix equations,
{\it SIAM J. Algebraic Discrete Methods}
{\bf 6} (1985),
406--412.
\end{enumerate}
\end{itemize}

\end{document}